\documentclass[a4paper]{amsart}
\usepackage{graphics}
\usepackage{color}
\usepackage{amssymb}
\usepackage[all]{xy}
\usepackage{amsmath}
\usepackage{stmaryrd}
\CompileMatrices
\newtheorem*{introthm}{Theorem}

\newtheorem{theorem}{Theorem}[section]
\newtheorem{lemma}[theorem]{Lemma}
\newtheorem{proposition}[theorem]{Proposition}
\newtheorem{corollary}[theorem]{Corollary}
\theoremstyle{definition}
\newtheorem{definition}[theorem]{Definition}
\newtheorem{example}[theorem]{Example}

\newtheorem{construction}[theorem]{Construction}

\newtheorem{remark}[theorem]{Remark}
\theoremstyle{remark}
\newtheorem*{introex}{Example}

\numberwithin{equation}{theorem}
\def\Chi{{\mathbb X}}

\def\div{{\rm div}}

\def\quot{/\!\!/}
\def\mal{\! \cdot \!}

\def\rq#1{\widehat{#1}}

\def\b#1{\overline{#1}}
\def\bangle#1{\langle #1 \rangle}

\def\KK{{\mathbb K}}
\def\TT{{\mathbb T}}
\def\ZZ{{\mathbb Z}}

\def\QQ{{\mathbb Q}}
\def\PP{{\mathbb P}}
\def\Of{{\mathcal{O}}}

\def\WDiv{\operatorname{WDiv}}

\def\Ample{{\rm Ample}}
\def\SAmple{{\rm SAmple}}
\def\rlv{{\rm rlv}}
\def\cov{{\rm cov}}
\def\Aut{\operatorname{Aut}}

\def\Cl{\operatorname{Cl}}

\def\Pic{\operatorname{Pic}}

\def\GL{{\rm GL}}

\def\Mat{{\rm Mat}}

\def\Spec{{\rm Spec}}

\def\cone{{\rm cone}}
\def\lin{{\rm lin}}

\def\topto#1{\stackrel{{\scriptscriptstyle #1}}{\longrightarrow}}
\def\faces{{\rm faces}}

\def\cov{{\rm cov}}
\def\SL{{\rm SL}}
\def\Sp{{\rm Sp}}
\def\Spin{{\rm Spin}}
\def\Ker{{\rm Ker}\,}
\def\Aut{{\rm Aut}}

\newcounter{itemnumber}

\begin{document}
\title[Embeddings with small boundary]
{On embeddings of homogeneous spaces\\ 
with small boundary}

\author[I.~Arzhantsev]%
{Ivan V. Arzhantsev} 
\thanks{Supported by DFG-Schwerpunkt 1094,
CRDF grant RM1-2543-MO-03, 
\\ 
\hspace*{\parindent}% 
RF President grant MK-1279.2004.1}
\address{Department of Higher Algebra, 
Faculty of Mechanics and Mathematics, 
Moscow State Lomonosov University,
Vorobievy Gory, GSP-2, Moscow, 119992, Russia}
\email{arjantse@mccme.ru}
\author[J.~Hausen]{J\"urgen Hausen} 
\address{Mathematisches Institut, Universit\"at T\"ubingen,
Auf der Morgenstelle 10, 72076 T\"ubingen, Germany}
\email{hausen@mail.mathematik.uni-tuebingen.de}
\begin{abstract}
We study equivariant embeddings with
small boundary of a given homogeneous space
$G/H$, where $G$ is a connected linear 
algebraic group 
with trivial Picard group and 
only trivial characters, 
and $H \subset G$ 
is an extension of a connected
Grosshans subgroup by a torus.
Under certain maximality conditions, 
like completeness, 
we obtain finiteness of the number 
of isomorphism classes of such embeddings, 
and we provide a combinatorial description 
of the embeddings and their morphisms.
The latter allows a systematic treatment
of examples and basic statements on the 
geometry of the equivariant 
embeddings of a given homogeneous space
$G/H$.
\end{abstract}

\maketitle

\section*{Introduction}

Homogeneous spaces $G/H$ and their 
equivariant (open) embeddings 
$G/H \subset X$ are of central interest 
in various fields of mathematics.
In the setting of algebraic geometry, 
there is a general approach to embeddings 
of homogeneous
spaces due to Luna and Vust~\cite{LuVu}.
However, this approach preferably applies 
to the case of small complexity, 
and even then, due to its generality, 
it is a deep and complicated theory, 
compare~\cite{Lu} and~\cite{Ti}.
In the present article, we concentrate
on a (rather) special class of 
$G/H$-embeddings, and for these we 
provide a simple alternative approach, 
based on combinatorial methods 
in Geometric Invariant Theory.

More precisely, let $G$ be a connected 
linear algebraic group
over an algebraically closed field $\KK$
of characteristic zero,
and assume that
$G$ has trivial Picard group
and admits only trivial characters
(e.g. is semisimple and simply connected).
We consider subgroups $H \subset G$,
which are ``Grosshans extensions'' in the 
sense that $H$ is an extension of a 
connected Grosshans subgroup $H_1 \subset G$ 
by some torus $T \subset G$; recall
from~\cite{Gr} that $H_1 \subset G$
is a Grosshans subgroup if and only if
$G/H_1$ is quasiaffine with a finitely 
generated algebra of global functions.

Given such a Grosshans extension
$H \subset G$, we investigate
``small'' equivariant embeddings $G/H \subset X$, 
where $X$ is a normal variety, and 
small means that the 
boundary $X \setminus (G/H)$ is small,
i.e., of codimension at least two in $X$.
Here are some simple examples.

\begin{introex}
For the special linear group 
$G := \SL(3,\KK)$, consider the 
connected Grosshans subgroup
\begin{eqnarray*}
H_1
& := &
\left\{
{\tiny
\left[
\begin{array}{ccc}
1 & 0 & a \\
0 & 1 & b \\
0 & 0 & 1
\end{array}
\right]
}; \; 
a,b \in \KK
\right\}.
\end{eqnarray*}
Then the following $G$-varieties
are small equivariant $G/H$-embeddings 
with a Gross\-hans extension $H \subset G$:
\begin{enumerate}
\item 
The product of projective spaces 
$\PP(\KK^3) \times \PP(\KK^3)$
with the diagonal $G$-action and the 
subgroup
\begin{eqnarray*}
H 
& := &
\left\{
{\tiny
\left[
\begin{array}{ccc}
t_1 & 0 & a \\
0 & t_2 & b \\
0 & 0 & t_1^{-1}t_2^{-1}
\end{array}
\right]
}; \; 
t_1,t_2 \in \KK^*, \; a,b \in \KK
\right\}
\end{eqnarray*}
\item 
The projective space $\PP(\KK^3 \oplus \KK^3)$
with $G$-action induced from the diagonal 
$G$-action on $\KK^3 \oplus \KK^3$ and the 
subgroup
\begin{eqnarray*}
H 
& := &
\left\{
{\tiny
\left[
\begin{array}{ccc}
t & 0 & a \\
0 & t & b \\
0 & 0 & t^{-2}
\end{array}
\right]
}; \; 
t \in \KK^*, \; a,b \in \KK
\right\}
\end{eqnarray*}
\end{enumerate}
\end{introex}

Note that in the first of these 
examples, the resulting
homogeneous space %$\SL(3)/H$ 
is spherical, i.e., of complexity zero, 
whereas in the second case it is 
of complexity one, and 
we have infinitely many $\SL(3)$-orbits.
However, also in the first setting, 
the construction may be generalized to
higher dimensions, see Proposition~\ref{xe1}, 
and then it produces $\SL(n)/H$-embeddings
of arbitrary high complexity, with 
$H \subset \SL(n)$ still being a 
Grosshans extension.

A first result of this paper is the 
following finiteness statement concerning 
isomorphism classes, see 
Theorem~\ref{cccc} for the full statement.

\begin{introthm}
Let $G$ be a connected 
linear algebraic group 
with trivial Picard group and
only trivial characters.
Then, for a fixed Grosshans extension $H \subset G$, 
there are only finitely many isomorphism
classes of complete small equivariant $G/H$-embeddings. 
%$X$ with small boundary $X \setminus (G/H)$.
\end{introthm}

Our main aim is to provide a combinatorial
description of the possible small equivariant 
$G/H$-embeddings $X$ and their morphisms
for a fixed Grosshans extension $H \subset G$.
This can be done under a maximality 
assumption: the variety $X$ should be $A_2$-maximal,
that means that any two points of $X$ 
admit a common 
affine neighbourhood in $X$, and
for every open embedding $X \subset X'$
into a variety $X'$ with the same 
property such that $X' \setminus X$ is 
of codimesion at 
least two in $X'$, we have $X=X'$.
Examples of $A_2$-maximal varieties are 
affine and projective ones, but there
exist definitely more of them.

Our approach is based on ideas of~\cite{BeHa1},
which we redevelop and enhance here
in a more geometric setting, 
see Sections~\ref{sec:1} and~\ref{sec:2}.
We consider the canonical action 
of the torus $T = H/H_1$ on the 
homogeneous space $G/H_1$. 
This action extends to the affine 
closure
\begin{eqnarray*}
Z 
& := & 
\Spec(\mathcal{O}(G)^{H_1}).
\end{eqnarray*}
The key observation is that 
every small equivariant 
$G/H$-em\-bedd\-ing~$X$ occurs as 
a good 
quotient space $U \quot T$ of a  
$T$-invariant open subset $U \subset Z$;
compare also~\cite{HK} for this
point of view.
Let us briefly see, what we obtain 
for the examples discussed before.

\begin{introex}[Continued]
For $G := \SL(3,\KK)$, and the Grosshans
subgroup $H_1 \subset G$ given as before,
we have
$$ 
Z 
\; = \;
\Spec(\mathcal{O}(G)^{H_1}) 
\; \cong \;
\KK^3 \oplus \KK^3.      
$$
Moreover, the open subsets $U \subset Z$ over 
(i) the product $\PP(\KK^3) \times \PP(\KK^3)$, 
and 
(ii) the projective space $\PP(\KK^3 \oplus \KK^3)$ 
are given by
$$ 
\text{(i)} \quad
U \; = \; \{(v_1,v_2); \; v_1 \ne 0 \ne v_2\},
\qquad
\text{(ii)}  \quad
U \; = \; \{(v_1,v_2); \; v_1 \ne 0 \text{ or }  v_2 \ne 0\}.
$$
\end{introex}

According to the key observation,
our task is now a matter of Geometric 
Invariant Theory: find an appropriate
description of the open $T$-invariant
sets $U \subset Z$ admitting a good quotient
$U \to U \quot T$ with an 
$A_2$-maximal quotient space.
Generalizing the description of projective 
quotients given in~\cite[Section~2]{BeHa2},
we present in Section~\ref{sec:1} a description  
in terms of ``orbit cones'' $\omega(z)$,
where $z \in Z$, which live in the 
rational character space $\Chi_\QQ(T)$.

Let us explain this.
For each $z \in Z$, define $\omega(z)$
to be the (convex, polyhedral) cone
generated by the weights $\chi \subset \Chi(T)$
admitting a semi-invariant 
$f \in \mathcal{O}(Z)$ with $f(z) \ne 0$.
To any open subset $U \subset Z$
we associate the collection 
$\Psi$ of the orbit cones $\omega(z)$,
where $T \mal z$ is closed in $U$.
It turns out that for the sets $U \subset Z$
with an $A_2$-maximal good quotient space
we obtain precisely the 2-maximal collections
$\Psi$, i.e., for any two cones of $\Psi$ 
their relative interiors overlap,
and $\Psi$ is maximal with this property.

Among the open subsets $U \subset Z$ with
an $A_2$-maximal good quotient space,
the small equivariant $G/H$-embeddings 
correspond to 
interior 2-maximal collections $\Psi$,
i.e., those containing the generic orbit 
cone.
Moreover, using the description of 
projective quotients in terms of the 
GIT-fan as provided in~\cite{BeHa2},
it is easy to figure out the 
projective small equivariant embeddings.
These observations are summarized 
in our second main result as follows,
see Theorem~\ref{embclass}.

\begin{introthm}
Let $G$ be a connected 
linear algebraic group 
with 
trivial Picard group and
only trivial characters,
and let $H \subset G$ be a Grosshans extension.
Then there is an equivalence of categories:
\begin{eqnarray*}
\left\{
\begin{array}{l}
\text{interior 2-maximal collections}
\\
\text{of orbit cones }
\end{array}
\right\}
& \to & 
\left\{
\begin{array}{l}
A_2 \text{-maximal small equivariant} 
\\
\text{} G/H \text{-embeddings}
\end{array}
\right\}
\end{eqnarray*}
If moreover $\mathcal{O}(G/H) = \KK$ holds, then we 
have in addition the following equivalence of categories:
\begin{eqnarray*}
\{
\text{interior GIT-cones}
\}
& \to & 
\left\{
\begin{array}{l}
\text{projective small equivariant} 
\\
G/H \text{-embeddings}
\end{array}
\right\}
\end{eqnarray*}
\end{introthm}

Strictly speaking, an equivalence of categories
needs a notion of morphisms on both sides.
On the left hand side, a morphism is a certain 
``face relation'', see Section~\ref{sec:2} for the
precise definition.
On the right hand side, we have, as usual, 
the equivariant base point preserving morphisms,
see Section~\ref{sec:3}. 

\begin{introex}[Continued]
Let $G := \SL(3,\KK)$ and 
$H_1 \subset G$ be as before.
Then, in the setting of~(i), the 
torus $T = H/H_1$ is of dimension 
two, and in $\Chi(T) = \ZZ^2$, 
the orbit cones of the $T$-action on 
$Z = \KK^3 \oplus \KK^3$ are 
the generic orbit cone
$$\omega(Z) \; = \; \cone(e_1,e_2) \; \subset \; \QQ^2 $$
and its faces. 
In particular, $\PP(\KK^3) \times \PP(\KK^3)$
is the only $A_2$-maximal small equivariant
$G/H$-embedding.
A little more variation takes place, 
if one considers a smaller torus extension 
of the Grosshans 
subgroup $H_1 \subset G$, for example:
\begin{eqnarray*}
H' 
& := &
\left\{
{\tiny
\left[
\begin{array}{ccc}
t & 0 & a \\
0 & t^{-1} & b \\
0 & 0 & 1
\end{array}
\right]
}; \; 
t \in \KK^*, \; a,b \in \KK
\right\}.
\end{eqnarray*}
The action of $T' = H' /H_1$ on $Z$ 
has besides the generic orbit cone 
$\QQ = \Chi_\QQ(T')$ the two  
rays $\QQ_{\le 0}$ and 
$\QQ_{\ge 0}$, and the 
zero cone $\{0\}$ as 
its orbit cones.
The latter three form the GIT-fan,
and the corresponding 
$G/H'$-embeddings are the only
$A_2$-maximal ones; 
see Example~\ref{hyperbolic}
for more details.
\end{introex}

As mentioned before, the construction 
of examples seen 
so far is put in Section~\ref{sec:4}
into a general framework. We provide 
several general constructions of 
spherical and non-spherical examples.
For simple groups $G$, we 
give in Section~\ref{sec:4} 
a classification 
of the small equivariant $G/H$-embeddings 
that additionally admit the structure of 
a toric variety, see 
Proposition~\ref{toricclass}.
Finally, we also construct some 
non-toric examples,
see Proposition~\ref{xe2}.

In the last Section, we 
study geometric properties of 
small equivariant $G/H$-embeddings.
The basic observation is that we may
apply the language of bunched rings
developed in~\cite{BeHa1}.
For example, existence 
of projective small embeddings 
with at most $\QQ$-factorial 
singularities immediately drops
out, see Corollary~\ref{existence}
and compare~\cite[Th\'{e}or\`{e}me~1]{BBII};
or one may produce 
non-projective complete 
small $G/H$-embeddings.
Moreover, we can easily construct 
examples of homogeneous spaces
$G/H$ admitting equivariant 
completions with small boundary
but no smooth ones, 
see Example~\ref{nosmoothemb}
and compare~\cite{BBK} for a discussion
of such phenomena.

\tableofcontents

%\newpage

\section{Quotients of affine torus actions}
\label{sec:1}

Given an algebraic variety with a reductive
group action, it is one of the basic 
tasks of Geometric Invariant Theory 
to describe all invariant open subsets
admitting a so called ``good quotient''.
The ``variation of quotients'' problem
is to understand the relations between 
these good quotients.
For reductive group actions on 
projective varieties, there is a 
satisfactory picture, concerning
good quotients that arise via Mumford's
method~\cite{Mu} from linearized
ample bundles, see~\cite{BrPr}, 
~\cite{DoHu}, \cite{Th}, and~\cite{Re}.

In~\cite[Section~2]{BeHa2}, the 
problem for quasiprojective 
quotient spaces of the action 
of a torus $T$ on a factorial affine 
variety $Z$ was considered, and 
an elementary construction of the
describing GIT-fan was given.
Here we ask more generally for open
subsets $U \subset Z$ that admit
a quotient space, which is 
embeddable into a toric variety.     
Our result generalizes a similar result
obtained in~\cite{BBSw} for linear 
representations of tori. 

We first fix our notation. 
By a lattice we mean a finitely generated
free abelian group.
For any lattice $K$, we denote by 
$K_\QQ := \QQ \otimes_\ZZ K$ the
associated rational vector space.
The word cone always stands for a 
convex, polyhedral cone in a rational
vector space.
For a cone $\sigma$, we denote its relative
interior by $\sigma^\circ$ and we write
$\tau \preceq \sigma$ if $\tau$ is a face
of $\sigma$.

In the sequel, 
$\KK$ is an algebraically
closed field, and, 
for the sake of rigorous
references, we suppose
$\KK$ to be characteristic zero
(though we expect our results 
to hold as well in positive 
characteristics).
Moreover, $K$ is a lattice, 
$T := \Spec(\KK[K])$ 
is the corresponding algebraic torus,
and we fix an $K$-graded integral 
affine $\KK$-algebra
\begin{eqnarray*}
R 
& = & 
\bigoplus_{u \in K} R_u.
\end{eqnarray*} 
Recall that the $K$-grading of $R$ defines 
a $T$-action on the corresponding 
affine variety $Z := \Spec(R)$ such that the 
homogeneous $f \in R_u$ are precisely the 
semi-invariants with respect to the character 
$\chi^u \colon T \to \KK^*$.

Now we recall some background around
good quotients; general references 
are~\cite{Se} and~\cite{BB2}.
A {\em good quotient\/} for a $T$-invariant
open set $U \subset Z$ is an affine, $T$-invariant 
morphism $\pi \colon U \to X$ such that the pullback 
map   
$\pi^* \colon \mathcal{O}_X \to \pi_*(\mathcal{O}_U)^T$
to the sheaf of invariants is an isomorphism.
If a $T$-invariant subset $U \subset Z$ admits a 
good quotient, then the quotient space is denoted
by $U \quot T$, and we will refer to $U \subset Z$ as a
{\em good $T$-set}.

It is a basic property of a good quotient 
$\pi \colon U \to U \quot T$ that
each of its fibres $\pi^{-1}(x)$ contains
precisely one closed $T$-orbit, and this 
orbit lies in the closure of any other
$T$-orbit of $\pi^{-1}(x)$.
From this one may derive the universal
property: any $T$-invariant morphism
$U \to Y$ factors uniquely through 
$U \to U \quot T$.  
This, by the way, justifies the notation 
$U \to U \quot T$.
One writes $U \to U/T$ for a good quotient, 
if it is {\em geometric}, i.e., its
(set-theoretical) fibres are precisely 
the $T$-orbits.

In the study of good $T$-sets the 
following concept is useful, 
compare~\cite{BB2}: 
an inclusion $U \subset U'$ of good 
$T$-sets in $Z$ is said to be 
{\em $T$-saturated\/} if    
$U$ is a full inverse image under 
the quotient map $U' \to U' \quot T$.
Due to the basic property of 
good quotients just mentioned, 
the set $U$ is $T$-saturated in 
$U'$ if and only if any closed 
$T$-orbit of $U$ is also closed 
in $U'$.

Let us define the good $T$-sets
$U \subset Z$ we are looking for.
First recall
that an {\em $A_2$-variety\/} is 
a variety $X$ such that any two 
$x,x' \in X$ admit a common 
affine neighbourhood in $X$.
For example, any quasiprojective 
variety is an $A_2$-variety.
It is shown in~\cite{Wl} that
the normal $A_2$-varieties are
precisely those admitting closed 
embeddings into toric varieties.

\begin{definition}[Compare~\cite{Sw2}]
We say that a good $T$-set $U \subset T$ 
is a {\em $(T,2)$-set\/} if the 
quotient space $U \quot T$ is an 
$A_2$-variety. 
By a {\em $(T,2)$-maximal\/}
subset of $Z$ we mean a $(T,2)$-set that 
does not occur as a proper $T$-saturated 
subset of some other $(T,2)$-set.
\end{definition}

Our aim is a combinatorial description
of all $(T,2)$-maximal subsets $U \subset Z$.
Let us introduce the necessary data.

As in~\cite{BeHa2}, we define the 
{\em orbit cone\/} associated to $z \in Z$
to be the (convex, polyhedral) cone 
$\omega(z) \subset K_\QQ$ generated
by all $u \in K$ that admit an element
$f \in R_u$ with $f(z) \ne 0$. 
There are only finitely many 
orbit cones, and we have 
$$ 
\omega(z)
\; = \; 
\omega(Z)
\; := \;
\cone(u \in K; \; R_u \ne 0)
$$
for all points $z$ of a nonempty open subset 
of $Z$.
Moreover, for any point $z \in Z$, the toric
variety $\Spec(\KK[\omega(z) \cap K])$
is the normalization of its $T$-orbit 
closure $C_Z(T \mal z)$ in $Z$.
In particular, we have a bijection
$$ 
\{T \text{-orbits in } 
C_Z(T \mal z)\}
\; \to \;
\faces(\omega(z)),
\qquad
T \mal z' \; \mapsto \; \omega(z').
$$

\begin{definition}
Let $\Omega(Z)$ denote the collection of
all orbit cones $\omega(z)$, where $z \in Z$.
\begin{enumerate}
\item
By a {\em 2-connected collection\/}
we mean a subcollection
$\Psi \subset \Omega(Z)$ 
such that
$\tau_1^\circ \cap \tau_2^\circ
\ne \emptyset$ 
holds for any two $\tau_1, \tau_2 \in \Psi$.
\item
By a {\em 2-maximal collection\/}, 
we mean a 2-connected collection, which is not 
a proper subcollection of any other 
2-connected collection.
\item
We say that a 2-connected collection $\Psi$ is a 
{\em face\/} of a 2-connected collection $\Psi'$ 
(written $\Psi \preceq \Psi')$, if for any 
$\omega' \in \Psi'$ there is an $\omega \in \Psi$
with $\omega \preceq \omega'$. 
\end{enumerate}
\end{definition}

Note that the 2-maximal collections 
form a partially ordered set with respect 
to the face relation defined above.
Here comes the link to the torus action.

\begin{definition}
To any collection $\Psi\subset \Omega(Z)$,
we associate a $T$-invariant subset
$U(\Psi) \subset Z$ as follows:
\begin{eqnarray*}
U(\Psi)
& := & 
\{ z \in Z; \; \omega_0 \preceq \omega(z) 
\text{ for some } \omega_0 \in \Psi \}.
\end{eqnarray*} 
\end{definition}

We state two basic observations.
The first one shows that the set $U(\Psi)$
associated to a 2-maximal collection
is a union of certain localizations, and thus 
is in particular open in $Z$.
The second one characterizes the closed $T$-orbits 
in $U(\Psi)$.

\begin{lemma}
\label{adjusted}
Let $\Psi$ be a 2-maximal collection
in $\Omega(Z)$.
Then any $z \in U(\Psi)$ admits an open
neighbourhood $U(z) \subset U(\Psi)$
such that for every $u \in \omega(z)^{\circ}$
there is an $n > 0$ and a homogeneous $f \in R_{nu}$
with $U(z) = Z_f$.
\end{lemma}

\begin{proof}
Choose homogeneous
$h_1, \ldots, h_r \in R$
such that $h_i(z) \ne 0$ holds
and the orbit cone
$\omega(z)$ is generated by 
$\deg(h_i)$, where $1 \le i \le r$.
For $w \in \ZZ_{>0}^r$, 
consider
$$ 
f^w
\; := \; 
h_1^{w_1} \ldots h_r^{w_r}
\; \in \; R.
$$
Then the sets $Z_{f^w}$ do not depend on 
the particular choice of  $w \in \ZZ_{>0}^r$.
Moreover for any $u \in \omega(z)^\circ$,
we find some $w$ with 
$$
\deg(f^w)
\; = \; 
w_1 \deg(h_1) +  \ldots + w_r \deg(h_r) 
\; \in \; 
\QQ_{> 0}u.
$$

In order to see that $U(z) = Z_{f^w}$ is 
as wanted, we still have to verify that
any $z' \in Z_{f^w}$ belongs to $U(\Psi)$.
By construction, we have 
$\omega(z) \subset \omega(z')$.
Consider $\omega_0 \in \Psi$ with 
$\omega_0 \preceq \omega(z)$. Then 
$\omega_0^\circ$ is contained in the 
relative interior of some face 
$\omega_0' \preceq \omega(z')$. 
By maximality of $\Psi$, 
we have $\omega_0' \in \Psi$,
and hence $z' \in U(\Psi)$.
\end{proof}

\begin{lemma}
\label{closedorbitchar}
Let $\Psi$ be a 2-connected collection in 
$\Omega(Z)$, and let $z \in U(\Psi)$.
Then the orbit $T \mal z$ is closed in 
$U(\Psi)$ if and only if 
$\omega(z) \in \Psi$ holds.
\end{lemma}

\begin{proof}
First let $T \mal z$ be closed in 
$U(\Psi)$.
By the definition of $U(\Psi)$, we have 
$\omega_0 \preceq \omega(z)$ for some 
$\omega_0 \in \Psi$.
Consider the closure $C_Z(T \mal z)$ of 
$T \mal z$ taken in $Z$, and choose
$z_0 \in C_Z(T \mal z)$ with 
$\omega(z_0) = \omega_0$.
Again by the definition of $U(\Psi)$,
we have $z_0 \in U(\Psi)$.
Since $T \mal z$ is closed in $U(\Psi)$, 
we obtain $z_0 \in T \mal z$,
and hence $\omega =\omega_0 \in \Psi$.

Now, let $\omega(z) \in \Psi$.
We have to show that any 
$z_0 \in C_Z(T \mal z) \cap U(\Psi)$
lies in $T \mal z$. 
Clearly, $z_0 \in C_Z(T \mal z)$ implies
$\omega(z_0) \preceq \omega(z)$.
By the definition of $U(\Psi)$, 
we have $\omega_0 \preceq \omega(z_0)$ 
for some $\omega_0 \in \Psi$.
Since $\Psi$ is a 2-connected collection,
we have 
$\omega_0^{\circ} \cap \omega(z)^\circ \ne \emptyset$.
Together with $\omega_0 \preceq \omega(z)$
this implies
$\omega_0 = \omega(z_0) = \omega(z)$,
and hence $z_0 \in T \mal z$. 
\end{proof}

A first major step towards the main
result of this section  
is to show that the 2-maximal 
collections define $(T,2)$-sets.

\begin{proposition}
\label{goodcoll2quot}
For any 2-maximal collection $\Psi$ 
in $\Omega(Z)$,
the associated $U(\Psi)$ is a $(T,2)$-set.
\end{proposition}

\begin{proof}
We regard $U(\Psi)$ as a union of 
sets $U(z)$ as provided in Lemma~\ref{adjusted}, 
where 
$z \in U(\Psi)$ runs through those points
that have a closed $T$-orbit in $U(\Psi)$;
according to Lemma~\ref{closedorbitchar}
these are precisely the points
$z \in U(\Psi)$ with $\omega(z) \in \Psi$.

First consider two such $z_1, z_2 \in U(\Psi)$.
Then we have $\omega(z_i) \in \Psi$, and
we can choose homogeneous 
$f_{1}, f_{2} \in R$ 
such that $\deg(f_1) = \deg(f_2)$ lies in 
$\omega(z_1)^\circ \cap \omega(z_2)^{\circ}$ 
and $U(z_i) = Z_{f_i}$ holds.
Thus, we obtain a commutative diagram
$$ 
\xymatrix{
{Z_{f_{1}}} \ar[d]^{\quot T} 
& {Z_{f_{1}f_{2}}} \ar[l]
                       \ar[r] 
                       \ar[d]^{\quot T} 
& {Z_{f_{2}}} \ar[d]^{\quot T} \\
{X_{f_{1}}} 
& {X_{f_{1}f_{2}}} \ar[l]
\ar[r]
& {X_{f_{2}}}
}
$$ 
where the upper horizontal maps are open embeddings,
the downwards maps are good quotients for the respective
affine $T$-varieties,
and the lower horizontal arrows indicate the induced
morphisms of the affine quotient spaces.

By the choice of $f_{1}$ and $ f_{2}$, 
the quotient $f_{2}/f_{1}$
is an invariant function on $Z_{f_{1}}$, and
the inclusion $Z_{f_{1}f_{2}} \subset Z_{f_{1}}$
is just the localization by $f_{2}/f_{1}$.
Since $f_{2}/f_{1}$ is invariant, the latter holds
as well for the quotient spaces; that means that
the map $X_{f_{1}f_{2}} \to X_{f_{1}}$ is localization
by $f_{2}/f_{1}$.

Now, cover $U(\Psi)$ by sets $U(z_i)$ with 
$T \mal z$ closed in $U(\Psi)$.
The preceding consideration allows gluing
of the maps $U(z_i) \to U(z_i) \quot T$ 
along
$U_{ij} \to U_{ij} \quot T$, where 
$U_{ij} := U(z_i) \cap U(z_j)$.
This gives a good quotient $U(\Psi) \to U(\Psi) \quot T$.
The quotient space is separated, because
we always have surjective multiplication maps
\begin{eqnarray*}
\mathcal{O}(Z_{f_{i}})^T \otimes \mathcal{O}(Z_{f_{j}})^T 
& \to &
\mathcal{O}(Z_{f_{i}f_{j}})^T.
\end{eqnarray*}

In order to see that $X = U(\Psi) \quot T$ is even 
an $A_2$-variety, consider $x_1,x_2 \in X$.
Then there are $f_i$ as above with $x_i \in X_{f_i}$.
The union $X_{f_1} \cup X_{f_2}$ is quasiprojective,
because, for example, the set 
$X_{f_1} \setminus X_{f_2}$ defines an ample divisor.
It follows that there is a common affine neighbourhood
of $x_1, x_2$ in $X_{f_1} \cup X_{f_2}$ and hence in 
$X$. 
\end{proof}

Now we discuss an inverse construction,
associating to any $T$-invariant open
subset a collection of orbit cones.

\begin{definition}
To any $T$-invariant open subset $U \subset Z$, 
we associate a set of orbit cones, namely
\begin{eqnarray*}
\Psi(U)
& := &  
\{\omega(z); \; z \in U \text{ with } 
T \mal z \text{ closed in } U
\}.
\end{eqnarray*}
\end{definition}

The following statement shows that,
when starting with a $(T,2)$-set, we 
obtain a 2-connected collection.
Its proof is the only place, where 
factoriality of the ring $R$ comes in;
in fact, it would even be sufficent 
to require that for any Weil divisor 
on $Z = \Spec(R)$ some positive 
multiple is principal.

\begin{proposition}
\label{T2goodcoll}
Suppose that $Z = \Spec(R)$ is factorial.
Then, for any $(T,2)$-set $U \subset Z$,
the associated set $\Psi(U)$ is a 
2-connected collection.  
\end{proposition}

\begin{proof}
By definition, the elements of $\Psi(U)$ 
are precisely the orbit cones $\omega(z)$, 
where $T \mal z$ is a closed subset of $U$.
We have to show that for any two cones
$\omega(z_i) \in \Psi(U)$, their relative interiors 
intersect nontrivially.

Consider the quotient $\pi \colon U \to U \quot T$,
and let $V \subset U \quot T$ be a common affine 
neighbourhood of $\pi(z_1)$ and $\pi(z_2)$.
Since $R$ is factorial, there is a homogeneous
function $f \in R$ vanishing precisely on 
the complement $Z \setminus \pi^{-1}(V)$.
It follows that the degree of $f$ 
lies in the relative
interior of both cones, $\omega(z_1)$
and $\omega(z_2)$.
\end{proof}

We are ready to formulate the main result of this 
section; it gives a complete description of the 
$(T,2)$-maximal sets $U \subset Z$, and describes
the possible inclusions of such sets.

\begin{theorem}
\label{A2maximal}
Let the algebraic torus 
$T = \Spec(\KK[K])$ act
on a factorial variety
$Z = \Spec(R)$.
Then we have mutually inverse 
bijections of finite sets:
\begin{eqnarray*}
\left\{
 \text{2-maximal collections in } \Omega(Z)
\right\}
& \longleftrightarrow &
\{(T,2) \text{-maximal } \text{subsets of } Z \}
\\
\Psi 
& \mapsto & 
U(\Psi)
\\
\Psi(U)
& \mapsfrom &
\Psi 
\end{eqnarray*}
These bijections are order-reversing maps of 
partially ordered sets in the sense 
that we always have 
\begin{eqnarray*}
\Psi \; \preceq \Psi'
& \Leftrightarrow &
U(\Psi) \; \supseteq U(\Psi').
\end{eqnarray*}
\end{theorem}

\begin{corollary}[\'Swi\c{e}cicka, \cite{Sw2}]
\label{T2finite}
The number of $(T,2)$-maximal
subets $U \subset Z$ is finite.
\end{corollary}

\begin{remark}
Theorem~\ref{A2maximal} shows that
parts of Bia\l ynicki-Birula's 
program~\cite{BB2}
can be carried out to obtain open 
subsets with good quotient for 
the action of a connected reductive 
group $G$ on $Z$: fix any maximal torus
$T \subset G$, determine the 
$(T,2)$-maximal open subset $U(\Psi)$
along the lines of Theorem~\ref{A2maximal},
and then~\cite[Theorem~1.1]{Ha1} 
provides good $G$-sets:
\begin{eqnarray*}
W(\Psi)
& = & 
\bigcap_{g \in G} g \mal U(\Psi).
\end{eqnarray*}
\end{remark}

\begin{proof}[Proof of Theorem~\ref{A2maximal}]
So far, we know from 
Proposition~\ref{goodcoll2quot}
that $U(\Psi)$ is a $(T,2)$-set,
and from Proposition~\ref{T2goodcoll} that 
$\Psi(U)$ is a 2-connected collection.
We begin with two auxiliary statements:

\medskip

\noindent
{\em Claim 1. } For any 2-maximal collection
$\Psi$ in $\Omega(Z)$
we have $\Psi(U(\Psi)) = \Psi$.

\medskip

Consider any $\omega \in \Psi(U(\Psi))$. 
By the definition of $\Psi(U(\Psi))$,
we have $\omega = \omega(z)$ for some 
$z \in U(\Psi)$
such that $T \mal z$ is closed in
$U(\Psi)$.
According to Lemma~\ref{closedorbitchar},
the latter implies $\omega \in \Psi$.
Conversely, let $\omega \in \Psi$.
Then we have $z \in U(\Psi)$,
for any $z \in Z$ with $\omega(z) = \omega$.
Moreover, Lemma~\ref{closedorbitchar}
tells us that $T \mal z$ is closed 
in $U(\Psi)$.
This implies $\omega \in \Psi(U(\Psi))$.

\medskip

\noindent
{\em Claim 2. }
Let $U \subset Z$ be a $(T,2)$-set,
and let $\Psi$ be any 2-maximal
collection in $\Omega(Z)$ 
with $\Psi(U) \subset \Psi$.
Then we have a $T$-saturated inclusion
$U \subset U(\Psi)$.

\medskip

First let us check that $U$ is 
in fact a subset of $ U(\Psi)$.
Given $z \in U$, we may choose 
$z_0 \in C_Z(T \mal z)$ such that 
$T \mal z_0$ is closed in $U$.
By definition of $\Psi(U)$, we
have $\omega(z_0) \in \Psi(U)$,
and hence $\omega(z_0) \in \Psi$.
Thus, $\omega(z_0) \preceq \omega(z)$
implies $z \in U(\Psi)$.

In order to see that the inclusion 
$U \subset U(\Psi)$ is $T$-saturated,
let $z \in U$ with $T \mal z$ closed in 
$U$.
We have to show that any
$z_0 \in C_Z(T \mal z)$
with $T \mal z_0$ closed in 
$U(\Psi)$
belongs to $T \mal z$.
On the one hand,
given such $z_0$, we obtain, 
using Claim~1:
$$ 
\omega(z_0) 
\; \in \; 
\Psi(U(\Psi))
\; = \;
\Psi.
$$
On the other hand, 
the definition of $\Psi(U)$
yields $\omega(z) \in \Psi$, 
and $z_0 \in C_Z(T \mal z)$
implies $\omega(z_0) \preceq \omega(z)$. 
Since $\Psi$ is a 2-connected collection,
the relative interiors of $\omega(z_0)$
and $\omega(z)$ intersect nontrivially,
and we obtain $\omega(z_0) = \omega(z)$.
This gives $z_0 \in T \mal z$. 

\medskip

Now we turn to the assertions of the Theorem.
First we show that the assignment 
$\Psi \mapsto U(\Psi)$ is well defined,
i.e., that $U(\Psi)$ is $(T,2)$-maximal.
Consider any $T$-saturated inclusion
$U(\Psi) \subset U$ with a $(T,2)$-set 
$U \subset Z$.
Using Claim~1, we obtain 
$$ 
\Psi
\; = \; 
\Psi(U(\Psi))
\; 
\subset 
\; 
\Psi(U).
$$ 
By maximality of $\Psi$, this implies
$\Psi = \Psi(U)$.
Thus, we obtain 
$U(\Psi) = U(\Psi(U))$.
By Claim~2, the latter set comprises
$U$, and thus, we see $U(\Psi) = U$.
In other words, $U(\Psi)$ is 
$(T,2)$-maximal.

Thus, we have a well-defined map 
$\Psi \to U(\Psi)$ from the 
2-maximal collections in $\Omega(Z)$
to the $(T,2)$-maximal subsets of $Z$.
According to Claim~1, this map 
is injective.
To see surjectivity, consider 
any $(T,2)$-maximal $U \subset Z$.
Choose a 2-maximal 
collection $\Psi$ with 
$\Psi(U) \subset \Psi$.
Claim~2 then shows $U = U(\Psi)$.
The fact that $\Psi \mapsto U(\Psi)$
and $U \mapsto \Psi(U)$ are inverse to 
each other is then obvious.

Finally, let us turn to the second statement
of the assertion.
The subset $U(\Psi')$ is contained in $U(\Psi)$ 
if and only if any closed $T$-orbit in $U(\Psi')$ 
is contained in $U(\Psi)$. 
By Lemma~\ref{closedorbitchar}, the points 
with closed $T$-orbit in $U(\Psi')$ are 
precisely the points $z \in Z$ with 
$\omega(z) \in \Psi'$.
By the definition of $U(\Psi)$,
such a point $z$ belongs to $U(\Psi)$
if and only if $\omega(z)$
has a face contained in $\Psi$.
\end{proof}

Finally, let us ask for
good $T$-sets $U \subset Z$ with a projective
quotient space $U \quot T$.
Clearly, such good $T$-sets are $(T,2)$-maximal, 
and we would like to see to which class of 
2-maximal collections they correspond. 

For this purposes, it is reasonable 
to assume that $R_0 = \KK$ holds, i.e., 
that there are only constant invariants.
Then the good $T$-sets $U \subset Z$ 
with a projective quotient space are 
precisely the sets of semistable points
of the $T$-linearizations of the trivial
line bundle over $Z$,
compare~\cite[Converse~1.12 and~1.13]{Mu}.

In~\cite{BeHa2}, the following description 
of the collection of the possible sets of 
semistable points is given: for any $u \in K$ 
define its {\em GIT-cone\/} to be the 
(convex, polyhedral) cone
\begin{eqnarray*}
\kappa(u)
& := & 
\bigcap_{u \in \omega(z)} \omega(z).
\end{eqnarray*} 
The main result of~\cite[Section~2]{BeHa2}
says that these GIT-cones form a fan,
and that they are in a (well defined) 
order reversing 
one-to-one correspondence to the possible 
sets of semistable points via
$$ 
\kappa \; \mapsto \; 
U(\kappa) := \bigcup_{f \in R_{nu}, n > 0} Z_f,
\quad \text{where } u \in \kappa^\circ.   
$$ 
So, given a GIT-cone $\kappa \subset K_\QQ$,
one may fix any $u \in K$ belonging 
to the relative interior of 
$\kappa^\circ \subset \kappa$,
and then $U(\kappa)$ is 
the union of all localizations
$Z_f$, where $f \in R$ is homogeneous of 
some degree $nu$ with $n > 0$.

\begin{proposition}
\label{pr8}
In the setting of Theorem~\ref{A2maximal},
suppose that $R_0 = \KK$ holds.
Then there is an order preserving 
injection of finite sets:
\begin{eqnarray*}
\left\{
GIT-cones
\right\}
& \longrightarrow &
\left\{
2 \text{-maximal collections in } \Omega(Z)
\right\}
\\
\kappa
& \mapsto & 
\Psi_\kappa \; := \; 
\{\omega \in \Omega(Z); \; \kappa^{\circ} \subset \omega^{\circ}\}.
\end{eqnarray*} 
The resulting open sets $U(\Psi_\kappa) \subset Z$ 
are precisely the good $T$-sets  in $Z$
that have a projective quotient space.
\end{proposition}

\begin{proof}
First recall that every good $T$-set 
$U \subset Z$ with a projective
quotient space $U \quot T$ is 
$(T,2)$-maximal.
Thus, our task is to show that for any 
of these $U \subset Z$ we have 
$\Psi(U) = \Psi_\kappa$,
with a unique GIT-cone $\kappa$.

Given a good $T$-set $U \subset Z$ with a 
projective quotient space, we know that
it is a set of semistable points, i.e.,
we have $U = U(\kappa)$ with a unique
GIT-cone $\kappa$.
Consider any $u \in \kappa^{\circ}$.
Then we have
$$
U
\; = \; 
U(\kappa)
\; = \; 
\bigcup_{f \in R_{nu}, n > 0} Z_f
\; = \; 
\{z \in Z; \; u \in \omega(z)\}.
$$
From the last description we infer 
that the closed $T$-orbits
of $U= U(\kappa)$ are precisely those 
$T \mal z \subset Z$, for which we have
$u \in \omega(z)^{\circ}$.
This implies
$\Psi(U) = \Psi_\kappa$,
and we are done.
\end{proof}

%\newpage

\section{Small $V$-embeddings}
\label{sec:2}

In this section, we prepare our study 
of equivariant embeddings of homogeneous 
spaces with small boundary.
We provide a framework for comparing 
varieties having a prescribed finitely
generated total coordinate ring.
But first, we recall the latter notion
and a little background,
compare~\cite{BeHa1} and~\cite{ElKuWa}.

Consider a normal variety $X$ with free finitely
generated divisor class group, and 
choose a subgroup $K \subset \WDiv(X)$
of the group of Weil divisors such that
the canonical map $K \to \Cl(X)$ is an
isomorphism.
The corresponding {\em total coordinate ring\/}
$\mathcal{R}_K(X)$ is defined as the algebra 
of global sections of a certain $K$-graded sheaf:
$$ 
\mathcal{R}_K(X)
\; := \;
\Gamma(X,\mathcal{R}_K), 
\qquad
\text{where} 
\quad
\mathcal{R}_K
\; := \;
\bigoplus_{D \in K} 
\mathcal{O}(D). 
$$
Note that for any homogeneous 
element $f \in \mathcal{R}_K(X)$
of degree $D \in K$ we obtain 
the sections over 
$X \setminus Z(f)$,
where $Z(f)$ is the suppport of 
$\div(f)+D$, as
\begin{eqnarray*}
\Gamma(X \setminus Z(f), \mathcal{R}_K)
& = &
\mathcal{R}_K(X)_f.
\end{eqnarray*}
This shows in particular that 
$\mathcal{R}_K$ is locally of finite
type, if $\mathcal{R}_K(X)$ is 
finitely generated.
Clearly, $\mathcal{R}_K$ is 
locally of finite type for any 
locally factorial $X$.

If $\mathcal{R}_K$ is locally of finite
type, then we may construct a 
(generalized) universal torsor: 
consider the relative spectrum 
$\rq{X} := \Spec_X(\mathcal{R}_K)$. 
The $K$-grading
of $\mathcal{R}_K$ defines an action of
the Neron-Severi torus $T := \Spec(\KK[K])$
on $\rq{X}$, and the canonical map 
$\rq{X} \to X$
is a good quotient for this action.

The variety $\rq{X}$ is quasiaffine.
If $\mathcal{R}_{K}(X)$ is finitely 
generated, then, setting 
$\b{X} := \Spec(\mathcal{R}_{K}(X))$,
and denoting by 
$F(X) \subset \mathcal{R}_K(X)$ the collection 
of homogeneous global sections such that 
$X \setminus Z(f)$ is affine, we obtain 
$\rq{X}$ as a union of localizations in 
$\b{X}$, compare~\cite[Proposition~3.1]{BeHa2}
and \cite[Proof of~4.2~a)]{BeHa1}:
$$
\rq{X}
\; = \;
\bigcup_{f \in F(X)} \b{X}_{f}
\; \subset \;
\b{X}.
$$

Now we are ready to begin our investigation
of varieties $X$ with prescribed finitely 
generated total coordinate ring $R$.
We start with an open subvariety 
$W \subset \Spec(R)$
that will serve as the universal torsor of
a ``common'' open subvariety $V$ of the 
varieties $X$.

%\begin{definition}
%For $x \in X$, we define its {\em divisorial cone\/}
%to be the cone in $\delta(x) \subset K_\QQ$ generated 
%by all divisors $D \in K$ admitting a section 
%$f \in \Gamma(X,\mathcal{O}(D))$ that does not 
%vanish near $x$, i.e.
%$$ 
%x 
%\; \not\in \; 
%Z(f)
%\; = \; 
%\Supp(\div(f) + D).
%$$ 
%\end{definition}

More precisely, let $W$ be a quasiaffine 
variety with trivial divisor class 
group $\Cl(W)$ such that $R := \mathcal{O}(W)$
satisfies $R^* = \KK^*$.
Moreover, suppose that there is a free action 
of an algebraic torus $T = \Spec(\KK[K])$ 
on $W$ admitting a geometric 
quotient $q \colon W \to V$.
This action determines 
a grading
\begin{eqnarray*}
R 
& = & 
\bigoplus_{u \in K}
R_u.
\end{eqnarray*}

We are now going to fix a group
of divisors  
$K_V \subset \WDiv(V)$
on the orbit space.
By freeness of the action on $W$,
we may fix a lattice basis 
$\{u_1, \ldots, u_k\}$ of $K$ such that 
for each $u_i$ there is a nonzero 
rational function $f_i \in Q(R)$, 
which is homogeneous of degree $u_i$.
We set
$$
f_u
\; := \; 
f_1^{n_1} \dots f_k^{n_k},
\quad 
\text{ for } 
u = n_1u_1 + \ldots + n_ku_k.
$$
Using once more freeness of the $T$-action on $W$, 
we may cover 
$W$ by $T$-invariant open subsets $W_\alpha \subset W$ 
such that for every $i = 1, \ldots, k$ there 
are invertible elements 
$\eta_{i,\alpha} \in \mathcal{O}(W_\alpha)_{u_i}$.
Similarly as before, we define 
$$
\eta_{u,\alpha}
\; := \; 
\eta_{1,\alpha}^{n_1} \dots \eta_{k,\alpha}^{n_k},
\quad 
\text{ for } 
u = n_1u_1 + \ldots + n_ku_k.
$$
This allows us to associate to any degree
$u \in K$ a Cartier divisors $D_u$ on the
orbit space $V = W/T$. 
Namely, denoting  $V_\alpha := q(W_\alpha)$, 
we define this divisor via local equations:
\begin{eqnarray*}
D_u \vert_{V_\alpha}
& := & 
\div\left(\frac{f_u}{\eta_{u,\alpha}} \right)
\end{eqnarray*}
We shall denote by $K_V \subset \WDiv(V)$ 
the group of divisors formed by the $D_u$,
where $u \in K$.
We list basic properties of this
construction, showing in particular that
$W$ is a universal torsor for the orbit
space $V = W/T$:

\begin{proposition}
\label{Wproperties}
The natural maps $K \to K_V$ and 
$K_V \to \Cl(V)$ are isomorphisms.
Moreover, we have a canonical isomorphism
of sheaves, given over open $V_0 \subset V$
by
$$ 
q_*\mathcal{O}_W \to \mathcal{R}_{K_V},
\qquad
\Gamma(V_0, q_* \mathcal{O}_W)_u \ni h
\; \mapsto  \;
\frac{h}{f_u} 
\in \Gamma(V_0,\mathcal{R}_{K_V})_u.
$$
This isomorphism of sheaves induces a 
canonical equivariant 
isomorphism $\rq{V} \to W$ of quasiaffine 
$T$-varieties. 
\end{proposition}

\begin{proof}
In order to see that $K \mapsto K_V$ is an 
isomorphism, we only have to care about injectivity.
This follows, for example, from $q^*(D_u) = \div(f_u)$ 
and $R^* = \KK^*$.
The fact that $K_V \to \Cl(V)$ is an isomorphism
is proven as in~\cite[Lemma~5.1]{BeHa1}, using standard 
arguments on computing the Picard group of 
a good quotient space provided in~\cite{KKV}.
The remaining statements are obvious.
\end{proof}

From now on, we assume that 
the $\KK$-algebra $R = \mathcal{O}(W)$ 
is finitely generated.
Then there is a canonical affine
closure $Z := \Spec(R)$.
The complement $Z \setminus W$
is of codimension at least two in 
$Z$, and $Z$ is a factorial variety.
Moreover,  the $T$-action on $W$ extends 
uniquely to a $T$-action on $Z$.

We consider {\em small $V$-embeddings}, 
i.e., open embeddings $\imath \colon V \to X$ into
a normal variety $X$ such that 
$X \setminus \imath(V)$ is of codimension at 
least two in $X$.
By a morphism of small $V$-embeddings
$\imath_1 \colon V \to X_1$ and
$\imath_2 \colon V \to X_2$, we mean 
a morphism $\varphi \colon X_1 \to X_2$ 
such that the following triangle is commutative:
$$ 
\xymatrix{
& 
V  \ar[dl]_{\imath_1} \ar [dr]^{\imath_2} 
& 
\\
X_1 \ar[rr]_{\varphi} & & X_2
}
$$

The following observation shows that 
small embeddings arise in a functorial 
way as quotient spaces of saturated 
extensions of the good $T$-set $W \subset Z$.

\begin{proposition}
\label{sat2small}
Let $U \subset U' \subset Z$ be
good $T$-sets, both containing
$W$ as a $T$-saturated subset.
Then we have:
\begin{enumerate}
\item
the induced map $V \to U \quot T$
of quotient spaces is a small $V$-embedding;
\item 
the induced map 
$U \quot T \to U' \quot T$
is a morphism of $V$-embeddings. 
\end{enumerate}
\end{proposition}

\begin{proof}
Since $W \subset U$ is a $T$-saturated
inclusion, the induced map of quotients
$V \to U \quot T$ is an open embedding. 
Moreover, this embedding must be
small, since $W \subset U$ is so.
The second statement is obvious.
\end{proof}

If $\imath \colon V \to X$ is a small
$V$-embedding, then any divisor $D \in K_V$ 
extends, by closing the support of 
$\imath_*(D)$, to a Weil divisor on $X$.
Denoting by $K_X \subset \WDiv(X)$ the group
of divisors obtained this way, we have 
canonical isomorphisms
$$ 
\Cl(V) 
\; \cong \; 
K_V 
\; \cong \; 
K_X
\; \cong \; 
\Cl(X) 
$$
The open embedding $\imath \colon V \to X$ induces 
a canonical isomorphism of graded sheaves 
$\imath^* \mathcal{R}_{K_X} \to \mathcal{R}_{K_V}$,
and hence we have an open embedding 
$\rq{\imath} \colon \rq{V} \to \rq{X}$. 
Moreover, the canonical isomorphism
$$ 
R \; \to \; \Gamma(X,\mathcal{R}_{K_X}),
\qquad
R_u \ni h 
\; \mapsto \;
(\imath^*)^{-1} \left(\frac{h}{f_u}\right)
\in
\Gamma(X,\mathcal{R}_{K_X})_u
$$
defines an isomorphism 
$\b{X} \to Z$.
The image of the open subset $\rq{X} \subset \b{X}$ 
is an open subset $W_X \subset Z$.
As we will see now,
the construction sending $X$ to $W_X$
is inverse to the one given in 
Proposition~\ref{sat2small}.

\begin{proposition}
\label{small2sat}
For every small embedding 
$\imath \colon V \to X$,
we have a $T$-saturated inclusion
$W \subset W_X$, and a commutative 
diagram
$$ 
\xymatrix@!0{
& Z & 
\\
W \ar[ur]^{\subset} \ar[rr]^{\subset} \ar[d]_{/T}
& &
{W_X} \ar[ul]_{\subset} \ar[d]^{\quot T} 
\\
V \ar[rr]_{\imath} 
& &
X
}
$$
Moreover, given  any 
$T$-saturated extension $W \subset U$ with
a good $T$-set $U \subset Z$, we have 
$W_{U \quot T} = U$.
\end{proposition}

\begin{proof}
Note that the isomorphism $\rq{V} \to W$ 
given in Proposition~\ref{Wproperties}
uniquely extends to an isomorphism $\b{V} \to Z$.
Thus, by construction, we have the following 
commutative diagram of $T$-equivariant maps:
$$ 
\xymatrix@!0{
& Z & 
\\
{\b{V}} \ar[ur]^{\cong} \ar[rr]^{\b{\imath}}
& &
{\b{X}} \ar[ul]_{\cong}  
\\
{\rq{V}} \ar[rr]_{\rq{\imath}} \ar[u]^{\subset} 
& &
{\rq{X}} \ar[u]_{\subset}
}
$$
This gives the commutative diagram in the 
assertion.
The fact that $W \subset W_X$ is 
$T$-saturated, follows from the fact that
this obviously holds for 
$\rq{\imath}(\rq{V}) \subset \rq{X}$.

Finally, let $W \subset U$ be a $T$-saturated 
extension. Set $X := U \quot T$, and 
consider the corresponding 
small embedding $V \to X$.
We have to show that $W_X = U$ holds.
Consider the commutative diagram
$$ 
\xymatrix@!0{
U \ar[ddr]_{q_U}
&
W \ar[l] \ar[r] \ar[d]
&
W_X \ar[ddl]^{q_X}
\\
& 
V \ar[d]
&
\\
& 
X 
& 
}
$$ 
Let $X_0 \subset X$ be any affine open 
subset. Then its boundary 
$D := X \setminus X_0$ is of pure 
codimension one in $X$.
Since $q_U$ as well as $q_X$ are affine,
and all sets of the upper row have small 
boundary in $Z$, we obtain
$$ 
q_U^{-1}(X_0)
\; = \; 
Z \setminus \b{q^{-1}(\imath^{-1}(D))}
\; = \; 
q_X^{-1}(X_0),
$$ 
where $\imath \colon V \to X$ 
denotes the embedding, and $q \colon W \to V$
is the geometric quotient fixed before.
The assertion now follows by covering
$X$ with affine open subsets $X_0 \subset X$.
\end{proof}

Next we note that morphisms of small
embeddings correspond to inclusions of 
$T$-saturated extensions of $W$ in $Z$.

\begin{proposition}
\label{Vembmorph}
Any morphism $\varphi \colon X_1 \to X_2$ 
of small $V$-embeddings $\imath_1 \colon V \to X_1$ 
and $\imath_2 \colon V \to X_2$ gives rise
to a commutative diagram 
$$
\xymatrix@!0{
& 
W \ar[dl] \ar [dr] \ar'[d][dd]
& 
\\
{W_{X_1}} \ar[rr] \ar[dd] & & {W_{X_2}} \ar[dd]
\\
& 
V  \ar[dl]_{\imath_1} \ar [dr]^{\imath_2} 
& 
\\
X_1 \ar[rr]_{\varphi} & & X_2
} 
$$
The map $W_{X_1} \to W_{X_2}$ is open 
inclusion. Moreover,   
$\varphi \colon X_1 \to X_2$ is an open 
embedding if and only if 
$W_{X_1} \subset W_{X_2}$
is $T$-saturated. 
\end{proposition}

\begin{proof}
The morphism $\varphi \colon X_1 \to X_2$ 
gives rise to a pullback homomorphism 
of sheaves
$\varphi^* \mathcal{R}_{X_2} \to \mathcal{R}_{X_1}$,
which in turn defines a commutative diagram
$$ 
\xymatrix{
& 
{\rq{V}} \ar[dr]^{\rq{\imath}_1} \ar[dl]_{\rq{\imath}_2} 
&
\\
{\rq{X}_1}
\ar[rr]_{\rq{\varphi}}
& 
&
{\rq{X}_1}
}
$$
Applying the canonical embeddings $\rq{V} \to Z$,
and $\rq{X}_i \to Z$, we obtain the commutative
diagram of the assertion.

The fact that the morphism $\varphi \colon X_1 \to X_2$ 
lifts to an inclusion
is due to the fact that $\imath_1$ and $\imath_2$ 
do so. 
Moreover, the fact that open embeddings 
$\varphi \colon X_1 \to X_2$
correspond 
to $T$-saturated inclusions $W_{X_1} \subset W_{X_2}$, 
follows from the observations that $W_{X_i} \to X_i$
is a good quotient and that open embeddings
$\varphi \colon X_1 \to X_2$ lift to $T$-saturated
open embeddings $\rq{\varphi} \colon \rq{X}_1 \to \rq{X}_2$.  
\end{proof}

The preceding three propositions may be summarized as
follows: on the one hand, we have the category of
saturated $W$-extensions, that means
good $T$-sets $U \subset Z$ containing $W$ as a 
$T$-saturated subset, together with inclusions as 
morphisms; on the other hand, we have the 
category of small $V$-embeddings. We showed:

\begin{corollary}
\label{smallembandsatext}
The assignments $U \mapsto U \quot T$ 
and $X \mapsto W_X$ define equivalences 
between the categories of saturated 
$W$-extensions and small $V$-embeddings, 
and they are essentially inverse to each
other.
\end{corollary}

As in~\cite{BeHa1}, we say that an
$A_2$-variety $X$ is {\em $A_2$-maximal\/}
if for any open embedding $X \to X'$ 
into an $A_2$-variety $X'$ such that
the complement $X' \setminus X$ is of codimension
at least two, we already have $X = X'$. 

\begin{corollary}
\label{A2maxsmallemb}
Let $\imath \colon V \to X$ be a small 
$V$-embedding into an $A_2$-variety $X$.
Then $X$ is $A_2$-maximal if and only if
$W_X \subset Z$ is $(T,2)$-maximal. 
\end{corollary}

As there exist complete normal varieties
which are not $A_2$-varieties, 
see~\cite[Example~6.4]{Sw2},
we would
like to have a similar statement 
comprising also the complete case.
Generalizing completeness, we introduce 
the following concept: 

\begin{definition}
W say that a variety $X$ is 
{\em 2-complete\/} if it admits no 
open embeddings $X \subset X'$ 
with $X' \setminus X$ nonempty 
and of codimension at least two.
\end{definition}

For small $V$-embeddings, also 2-completeness
may be expressed in terms of quotients. 
Recall from~\cite{BB2} that a good 
$T$-set $U \subset Z$ is said to be 
{\em T-maximal\/} if it is maximal with 
respect to $T$-saturated 
inclusion in the collection of 
all good $T$-sets of $Z$.

\begin{corollary}
\label{2completechar}
Let $\imath \colon V \to X$ be a small 
$V$-embedding. 
Then $X$ is 2-complete if and only if
$W_X \subset Z$ is $T$-maximal. 
\end{corollary}

%\newpage

\section{Small equivariant embeddings}
\label{sec:3}

In this section, we present the main 
results of the paper.
We study small equivariant
embeddings $G/H \subset X$,
where $H \subset G$ is a  
``Grosshans extension''; 
the precise definition is given below.
Among other things, we obtain finiteness
of the numbers  of isomorphism classes 
of 2-complete and $A_2$-maximal 
small equivariant $G/H$-embeddings,
and we give a combinatorial description
in the latter case.

Let us fix the setup.
Throughout this section, we denote 
by $G$ a connected linear algebraic 
group having a trivial character group $\Chi(G)$,
and trivial Picard group $\Pic(G)$.
For example, $G$ might be any 
connected simply connected semisimple group,
like the special linear group
$\SL(n,\KK)$.

Let $H \subset G$ be a closed subgroup.
As usual, we mean by an {\em equivariant embedding\/}
of the homogeneous space $G/H$ an irreducible, 
normal $G$-variety $X$ together with a base point 
$x_0 \in X$ such that $H$ equals the isotropy 
group $G_{x_0}$ of $x_0 \in X$ and the morphism 
$$ G/H \to X, \qquad gH \mapsto g \mal x_0 $$
is a ($G$-equivariant) open embedding.
A {\em morphism\/} of two equivariant embeddings
$X$ and $X'$ of $G/H$ is a $G$-equivariant morphism 
$X \to X'$ sending the base point $x_0 \in X$ to 
the base point $x'_0 \in X'$.
Note that if a morphism of $G/H$-embeddings exists, 
then it is unique.

\begin{definition}
By a {\em small equivariant $G/H$-embedding}
we mean an equivariant $G/H$-embedding 
with a normal variety $X$ such that the complement 
$X \setminus G \mal x_0$ is of codimension
at least two in $X$. 
\end{definition}

Let us recall from~\cite{Gr} the necessary 
concepts from algebraic group theory.
The subgroup $H \subset G$ is
said to be {\em observable\/} if $G/H$ 
is a quasiaffine variety.
Moreover, $H \subset G$ is called a {\em Grosshans 
subgroup\/} if it is observable and the 
algebra of global functions 
$\mathcal{O}(G/H) = \mathcal{O}(G)^H$ 
is finitely generated.

\begin{remark}
\label{Grosshansexamples}
In each of the following cases,
the subgroup $H \subset G$ is
a Grosshans subgroup:
\begin{itemize}
\item 
$G/H$ is quasiaffine and spherical or 
of complexity one, see~\cite{Kn};
\item
$H$ is the unipotent radical of a parabolic
subgroup of $G$, see~\cite{Gr};
\item
$H$ is the generic stabilizer of a 
factorial affine $G$-variety, see~\cite{Gr}.
\end{itemize}
\end{remark}

The property of being a Grosshans subgroup 
can (tautologically) be characterized in terms of 
small embeddings; more precisely, 
we observe the following.

\begin{remark}
\label{Grosshanschar}
The subgroup $H \subset G$ is Grosshans
if and only if there is a
small embedding $G/H \to X$ into a 
(normal) affine variety $X$. 
In this case, $X$ is the spectrum 
of $\mathcal{O}(G/H)$.  
\end{remark}

We are ready to introduce 
the notion of a ``Grosshans extension''.
Let $K := \Chi(H)$ denote the 
character group of the subgroup 
$H \subset G$.
Then we have an associated 
diagonalizable group 
$T := \Spec(\KK[K])$,
a canonical epimorphism
$\pi \colon H \to T$,
and we may consider its kernel:
$$ 
H_1 
\; := \;
\ker(\pi)
\; = \; 
\bigcap_{\chi \in K}
\ker(\chi)
\; \subset \;
H. 
$$ 

\begin{definition}
We say that $H \subset G$ is 
{\em a Grosshans extension},
if $H$ is connected,
and $H_1 \subset G$ is a 
Grosshans subgroup.
\end{definition}

We now present our results; the proofs
are given at the end of the section.
The first observation is a 
characterization of Grosshans extensions
in the spirit of Remark~\ref{Grosshanschar}.

\begin{proposition}
\label{grossextchar}
A connected closed subgroup $H \subset G$
is a Grosshans extension if and only if
there is a small embedding 
$G/H \to X$ into a
normal variety $X$ with finitely generated
free divisor class group and 
finitely generated total coordinate 
ring.
\end{proposition}

The next statement shows that, 
in many cases,
for a Grosshans extension $H \subset G$,
small embeddings $G/H \to X$ are 
automatically equivariant.

\begin{proposition}
\label{small2equivar}
Let $H \subset G$
be a Grosshans extension,
and let 
$\imath \colon G/H \to X$  
be a small embedding.
If $X$ is 2-complete or 
$A_2$-maximal, 
then the canonical $G$-action 
on $\imath(G/H) \subset X$ 
extends to the whole $X$. 
\end{proposition}

We come to the main results. 
The first one is the following 
finiteness statement on small equivariant 
embeddings.

\begin{theorem}
\label{cccc}
Let $G$ be a connected 
linear algebraic group 
with 
trivial Picard group and
only trivial characters,
and let $H \subset G$ be a Grosshans 
extension.
\begin{enumerate}
\item 
The number of isomorphism classes of 
2-complete small equivariant 
$G/H$-embeddings is finite. 
\item
The number of isomorphism classes of 
$A_2$-maximal small equivariant 
$G/H$-embeddings is finite. 
\end{enumerate}
\end{theorem}

As a direct application, we note the 
following statement on the group of
$G$-equivariant automorphisms:

\begin{corollary}
\label{norm2aut}
In the setting of~\ref{cccc},
let $N_G(H)$ be the normalizer 
of $H$ in $G$, 
and $N^0 \subset N_G(H)/H$ the 
unit component. 
Then, for any 2-complete or $A_2$-maximal 
small equivariant $G/H$-embedding $X$,
the group  $\Aut_G(X)$ of 
$G$-equivariant automorphisms
contains $N^0$ as its unit component.
 \end{corollary}

Now, fix a Grosshans
extension $H \subset G$.
Note that $K = \Chi(H)$ is a lattice,
and thus $T = \Spec(\KK[K])$ is a torus.
Consider the subgroup $H_1 \subset H$ as 
defined before.
Then $W := G/H_1$ is a quasiaffine 
variety, and it comes with canonical
actions of $G$ and $T \cong H/H_1$,
given by
$$ 
g \mal (g'H_1) 
\; := \; 
gg'H_1,
\qquad
(hH_1) \mal g'H_1
\; := \; 
g'h^{-1}H_1. 
$$
The algebra $R := \mathcal{O}(W)$ is finitely 
generated, and we have a canonical 
affine closure $Z := \Spec(R)$.
The actions of $G$ and $T$ both extend 
to $Z$, and, obviously, they commute.
The $T$-action on $Z$ defines a 
grading of the algebra of functions:
\begin{eqnarray*}
R
& := &
\bigoplus_{u \in K} R_{u}.
\end{eqnarray*}
Having in mind this grading,
we may speak, as in the first section,
about the collection $\Omega(Z)$
of all $T$-orbit cones 
$\omega(z) \subset K_\QQ$, and 
also about the GIT-cones 
$\kappa \subset K_\QQ$. 
We will work with the following
notions. 

\begin{definition}
We say that a subset $\Psi \subset \Omega(Z)$
is an {\em interior collection\/}
if it contains the weight cone $\omega(Z)$.
By  an {\em interior GIT-cone\/} we mean 
a GIT-cone $\kappa \subset K_\QQ$ with
$\kappa^\circ \subset \omega(Z)^\circ$.
\end{definition}

Our main result is a description of 
the category of $A_2$-maximal equivariant small 
$G/H$-embeddings of a given Grosshans extension
$H \subset G$ in terms of interior 2-maximal
collections $\Psi \subset \Omega(Z)$.
In the case $\mathcal{O}(G/H) = \KK$, 
it comprises also a description of all 
projective equivariant small 
$G/H$-embeddings.

Together with the face relations ``$\preceq$''
as morphisms, 
the interior 2-maximal collections 
$\Psi \subset \Omega(Z)$ form a category.
Moreover, we may associate to any $\Psi$ 
the variety $U(\Psi) \quot T$, and this assignment
is functorial: if we have $\Psi \preceq \Psi'$, 
then there is an induced morphism of the 
quotient spaces 
\begin{eqnarray*}
U(\Psi') \quot T 
& \to & 
U(\Psi) \quot  T.
\end{eqnarray*}
As we will see, the set $U(\Psi) \subset Z$ is 
$G$-invariant, and thus the $G$-action descends
to the quotient space $U(\Psi) \quot T$.
Moreover this space comes with a canonical 
base point, namely $\pi(e_GH_1)$, where
$\pi \colon U(\Psi) \to U(\Psi) \quot T$ 
denotes the quotient map.

\begin{theorem}
\label{embclass}
Let $G$ be a connected  
linear algebraic group 
with 
trivial Picard group and
only trivial characters, 
and let $H \subset G$ be a Grosshans extension.
Then we have a contravariant equivalence of 
categories:
\begin{eqnarray*}
\{
\text{interior 2-maximal collections}
\}
& \to & 
\left\{
\begin{array}{l}
A_2 \text{-maximal small equivariant} 
\\
\text{} G/H \text{-embeddings}
\end{array}
\right\}
\\
\\
\Psi
& \mapsto &
U(\Psi) \quot T.
\end{eqnarray*}
If moreover $\mathcal{O}(G/H) = \KK$ holds, then we 
have in addition a contravariant equivalence of categories:
\begin{eqnarray*}
\{
\text{interior GIT-cones}
\}
& \to & 
\left\{
\begin{array}{l}
\text{projective small equivariant} 
\\
G/H \text{-embeddings}
\end{array}
\right\}
\\
\\
\kappa
& \mapsto &
U(\kappa) \quot T.
\end{eqnarray*}
\end{theorem}

The condition $\mathcal{O}(G/H) = \KK$ has been 
studied by several authors; for example, 
various characterizations of this property 
and concrete examples can be found 
in~\cite{BBI}, \cite{BBII}, 
and~\cite[Section~23\,B]{Gr}.

\begin{remark}
If we drop the assumption $\mathcal{O}(G/H) = \KK$
in the second part of Theorem~\ref{embclass},
then the interior GIT-cones correspond to those
$A_2$-maximal small equivariant $G/H$-embeddings,
which are in addition quasiprojective.
\end{remark}

Using Proposition~\ref{pr8},
we observe that in the case of 
a small character group $\Chi(H)$
and $\mathcal{O}(G/H) = \KK$, 
every $A_2$-maximal small 
$G/H$-em\-beddings is projective
More precisely, we obtain the following.

%\begin{corollary}
%A Grosshans extension $H \subset G$
%with $\Chi(H) \cong \ZZ$ and $\mathcal{O}(G/H) = \KK$, admits 
%up to isomorphism exactly one 
%projective small equivariant $G/H$-embedding.
%\end{corollary}

\begin{corollary}
\label{Z2proj}
For a Grosshans extension $H \subset G$
with $\Chi(H)$ of rank at most two
and $\mathcal{O}(G/H) = \KK$,
every $A_2$-maximal small equivariant
$G/H$-embedding is projective.
\end{corollary}

Let us begin to prove the results.
A first basic step is the 
following group theoretical observation.

\begin{proposition}
\label{H1props}
Let $G$ be any linear algebraic group,
and let $H \subset G$ be a closed subgroup.
Then the following subgroup is observable:
\begin{eqnarray*}
H_1
& = &
\bigcap_{\chi \in \Chi(H)} \ker(\chi).
\end{eqnarray*}
Moreover, if the subgroup 
$H \subset G$ is connected,
then the above group $H_1$ has only 
trivial characters.
\end{proposition}

\begin{proof}
By Chevalley's theorem, 
there exist
a rational finite-dimensional 
$G$-module $V$,
and a non-zero vector $v\in V$
such that $H$ is the stabilizer
of the line $\KK v$.
Let $\chi_0 \in\Chi(H)$ such that
$h \mal v = \chi_0(h) v$ for any 
$h\in H$. 
Then $H_0 := \Ker(\chi_0)$ 
equals the isotropy subgroup $G_v$ 
of $v$, 
and thus is observable in $G$,
see~\cite[Theorem~2.1]{Gr}.
Note that we have
$$ 
H_1
\;  = \; 
\bigcap_{\chi \in \Chi(H)} \ker(\chi \vert_{H_0})
\; \subset \; 
H_0.
$$
Since the intersection of observable sugroups
is again observable, it suffices to show
that each $H_\chi := \ker(\chi \vert_{H_0})$ 
is observable.
For this, use again~\cite[Theorem~2.1]{Gr} 
to realize the one dimensional $H_0$-module
given by $\chi \vert_{H_0}$ as an $H_0$-submodule
$\KK v_\chi$ of a $G$-module $V_\chi$.
Then $H_\chi$ is the isotropy group
of $(v,v_\chi)$ in the $G$-module 
$V \oplus V_\chi$, and hence it is 
observable.   

To see the second assertion,
let $H=LR_u(H)$ be the Levi decomposition, 
where $R_u(H)$ is the unipotent radical 
and $L$ is connected reductive,
compare~\cite[Sec.~6.4]{OV}. 
The subgroup $L$ locally splits into 
the direct product $L=T^cL^s$, 
where $T^c$ is the central torus and
$L^s$ is a semisimple subgroup coinciding 
with the commutator subgroup of $L$, 
compare~\cite[Sec.~4.1.3]{OV}.
Clearly, $R_u(H)$ and $L^s$ are contained in $H_1$. 
On the other hand, 
$H/(L^sR_u(H))$ is isomorphic to 
$T^c/(T^c\cap L^s)$, and hence 
is a torus.
This implies $H_1=L^sR_u(H)$, 
which proves $\Chi(H_1)=0$. 
\end{proof} 

The next step is to ensure that we are in the 
setup of the preceding two sections.
Let $H \subset G$ be any connected subgroup,
set $W := G/H_1$ and $V := G/H$.

\begin{lemma}
\label{setup}
The variety $W$ is quasiaffine, factorial,
and satisfies $\mathcal{O}^*(W) = \KK^*$.
Moreover, the $T$-action on $W$ is free, 
and the canonical map $W \to V$ 
is a geometric quotient for this action.
In particular, $W \to V$ is a universal
torsor. 
\end{lemma}

\begin{proof}
Proposition~\ref{H1props} tells us that
$H_1 \subset G$ is a Grosshans subgroup,
and hence $W = G/H_1$ is quasiaffine.
It is obvious that $T \cong H/H_1$ acts freely
on $W = G/H_1$. 
Thus, the canonical map $W \to V$ must be a 
geometric quotient for the $T$-action on $W$.
Moreover $\mathcal{O}^*(W) = \KK^*$ follows from 
$\mathcal{O}^*(G) = \KK^*$, which in turn
is due to $\Chi(G) = 0$, see for 
example~\cite[Prop.~1.2]{KKV}.

It remains to show that $W$ is factorial, 
i.e., has trivial divisor class group $\Cl(W)$.
Since $W$ is smooth, we have $\Cl(W) = \Pic(W)$.
The latter group occurs in the 
exact sequence
$$ 
\xymatrix{
\Chi(G) \ar[r]
&
\Chi(H_1) \ar[r]
&
\Pic(W) \ar[r]
&
\Pic(G),
}
$$ 
see, for example~\cite[Prop.~3.2]{KKV}.
We assumed $\Chi(G) = 0$ and $\Pic(G) = 0$.
Thus $\Pic(W) = 0$ follows from $\Chi(H_1) = 0$.
Proposition~\ref{Wproperties} then tells us
that $W \to V$ is a universal torsor. 
\end{proof}

\begin{proof}[Proof of Proposition~\ref{grossextchar}]
The ``only if'' part follows directly from
Lemma~\ref{setup}: it tells us that $G/H_1 \to G/H$ 
is a universal torsor and hence the trivial embedding 
$G/H \to G/H$ is as wanted.
 
Conversely, if there is a small embedding 
$G/H \to X$
as in the assertion, then $G/H$ has a free
finitely generated divisor class group and a 
finitely generated total coordinate ring 
$\mathcal{R}(G/H)$.
Again Lemma~\ref{setup} shows that 
$G/H_1 \to G/H$ is a universal torsor.
Thus $\mathcal{O}(G/H_1) \cong \mathcal{R}(G/H)$ 
holds. In particular, this algebra is finitely
generated, which means that  
$H_1$ is a Grosshans subgroup. 
\end{proof}

From now on, $H \subset G$ is a Grosshans 
extension, and, thus $R = \mathcal{O}(W)$
is finitely generated.
Note that then $Z = \Spec(R)$ contains
$W$ as an open subset with small complement 
$Z \setminus W$.
In particular, we have $\Cl(Z) = 0$, which 
means that $Z$ is factorial.

\begin{lemma}
\label{Ginvariant}
Let $\Psi \subset \Omega(Z)$ be a 2-maximal
collection of orbit cones.
Then the associated set $U(\Psi) \subset Z$ 
is $G$-invariant. 
\end{lemma}

\begin{proof}
Since the actions of $T$ and $G$ on $Z$ commute,
we see that for any $z \in Z$ and any $g \in G$,
we have $\omega(z) = \omega(g \mal z)$.
The assertion thus follows from the definition
of $U(\Psi)$. 
\end{proof}

\begin{proof}[Proof of Proposition~\ref{small2equivar}]
Let us first consider the $A_2$-maximal case.
According to Corollary~\ref{A2maxsmallemb},
the small embedding $\imath \colon G/H \to X$, defines a 
$(T,2)$-maximal open subset $W_X \subset Z$
such that everything fits into a
commutative diagram
$$ 
\xymatrix{
W \ar[r]^{\subset} \ar[d]_{/T}
& 
{W_X} \ar[d]^{\quot T}
\\
G/H \ar[r]_\imath
&
X
}
$$
where $W \subset W_X$ is a $T$-saturated 
inclusion.
By Theorem~\ref{A2maximal},
we have $W_X = U(\Psi)$ 
with some 2-maximal collection
$\Psi \subset \Omega(Z)$. 
By Lemma~\ref{Ginvariant},
the set $W_X$ is $G$-invariant.
Thus, the action of $G$
on $W_X$ descends to the
desired $G$-action on $X$.

If $X$ is 2-complete, then
the arguments are similar. 
Again, by Proposition~\ref{small2sat}, 
we have a $T$-saturated inclusion 
$W \subset W_X$, where $W_X \subset Z$ 
is a good $T$-set.
Moreover, Corollary~\ref{2completechar}
tells us that $W_X$ is $T$-maximal.
Thus,~\cite[Corollary~2.3]{Sw1} yields 
that $W_X$ is $G$-invariant, and, 
again, the $G$-action 
descends to the desired action on 
the variety~$X$. 
\end{proof}

\begin{proof}[Proof of Theorem~\ref{cccc}]
By Proposition~\ref{small2equivar}, 
the category of 2-complete small equivariant
$G/H$-embeddings
as well as the category of 
$A_2$-maximal small equivariant
$G/H$-embeddings
are full subcategories of the category of
small (not necessarily equivariant) 
$V$-embeddings, where $V = G/H$.

Thus, according to 
Corollaries~\ref{smallembandsatext},
\ref{A2maxsmallemb}, and~\ref{2completechar} 
we only need to know that the collection of 
$T$-maximal open subsets $U \subset Z$
and the collection of $(T,2)$-maximal open 
subsets $U' \subset Z$ are finite.
In the first case this follows from the 
main result of~\cite{BB1}, in the second
case this is Corollary~\ref{T2finite}.
\end{proof}

\begin{proof}[Proof of Corollary~\ref{norm2aut}]
First note that $N_G(H)/H$ may be identified 
with the group $\Aut_G(G/H)$;
in fact, 
$N_G(H)/H$ acts on $G/H$ via 
\begin{eqnarray*}
nH \mal gH
& := & 
gn^{-1}H. 
\end{eqnarray*}
Consequently, for any $G$-variety
$X$ with an open $G$-orbit isomorphic 
to $G/H$, the group $\Aut_G(X)$ 
is a subgroup of $N_G(H)/H$.

Moreover, the group $N_G(H)/H$
acts on the set of isomorphism 
classes of $G/H$-embeddings
via
\begin{eqnarray*}
nH \mal (X,x_0)
& := &
(X,n^{-1} \mal x_0).
\end{eqnarray*}
Two pairs $(X,x_0)$ and $(X,n^{-1} \mal x_0)$ 
are isomorphic as $G/H$-embeddings if and only if 
$nH\in\Aut_G(X)$ holds.

For 2-complete and as well 
$A_2$-maximal equivariant
 $G/H$-embeddings $X$,
Theorem~\ref{cccc} tells us 
that the respective numbers 
of isomorphism classes 
 are finite.
Hence, for a given $X$, the group
$\Aut_G(X)$ acts on the set
of $H$-fixed points in $G/H$ 
with finitely many orbits. 
This action is precisely the action 
of $\Aut_G(X)$ on the group $N_G(H)/H$
by right multiplication, and thus,
by dimension reasons, $\Aut_G(X)$
contains $N^0$ as its unit component.
\end{proof}

\begin{proof}[Proof of Theorem~\ref{embclass}]
Let us first check that the assignment 
is well defined.
By Lemma~\ref{Ginvariant}, the sets 
$U(\Psi) \subset Z$ defined by 2-maximal 
collections $\Psi \subset \Omega(Z)$ are 
$G$-invariant.
Moreover, any interior 2-maximal  
$\Psi \subset \Omega(Z)$ contains 
the generic orbit cone $\omega(e_G H_1)$.
Consequently, $W \subset U(\Psi)$ holds, and,
by Lemma~\ref{closedorbitchar}, this is a 
$T$-saturated inclusion.
Proposition~\ref{sat2small} thus
provides a commutative diagram
$$ 
\xymatrix{
W \ar[r]^{\subset} \ar[d]_{/T}
& 
U(\Psi) \ar[d]^{\quot T}_{\pi}
\\
G/H \ar[r]
&
X
}
$$
where $X := U(\Psi) \quot T$, 
and the induced map of quotient spaces
is an open embedding.
The $G$-action on $U(\Psi)$ 
descends to an action on the 
quotient variety $X$, 
making it into a small equivariant 
$G/H$-embedding with base point
$\pi(e_G H_1)$.
By Corollary~\ref{A2maxsmallemb}
the variety $X$ is $A_2$-maximal.
Hence, the assignment 
$\Psi \mapsto U(\Psi) \quot T$
is well defined.

According to Proposition~\ref{small2equivar},
the category of $A_2$-maximal small 
equivariant $G/H$-embeddings 
and that of $A_2$-maximal small 
$G/H$-embeddings are 
isomorphic
via sending the $G$-variety $(X,x_0)$ 
to the embedding $G/H \to X$,
$gH \mapsto g \mal x_0$.
Thus, Theorem~\ref{A2maximal}, 
together with Corollaries~\ref{smallembandsatext}
and~\ref{A2maxsmallemb} 
shows that $\Psi \mapsto U(\Psi) \quot T$
defines a (contravariant)
fully faithful functor. 

Similarly as 
in the proof of Proposition~\ref{small2equivar}, 
we show now that our functor 
is essentially surjective.
Let $X$ be any
$A_2$-maximal small equivariant 
$G/H$-embedding.
Then, according to  
Corollary~\ref{A2maxsmallemb},
we have a commutative diagram
of $T$-equivariant maps with a 
$(T,2)$-maximal subset
$W_X \subset Z$ and a $T$-saturated
inclusion $W \subset W_X$:
$$ 
\xymatrix{
W \ar[r]^{\subset} \ar[d]_{/T}
& 
{W_X} \ar[d]^{\quot T}
\\
V \ar[r]
&
X
}
$$
By Theorem~\ref{A2maximal},
we have $W_X = U(\Psi)$ for a
2-maximal collection of orbit cones
$\Psi \subset \Omega(Z)$.
Moreover, by Lemma~\ref{closedorbitchar}, 
the generic orbit cone $\omega(e_GH_1)$
belongs to~$\Psi$.
Hence $\Psi$ is an interior
2-maximal collection.
Clearly, there is an induced 
isomorphism $X \to U(\Psi) \quot T$
of small equivariant $G/H$-embeddings.

The statement concerning the projective 
case is an immediate consequence of 
Proposition~\ref{pr8} and the 
$A_2$-maximal case, which
we just settled. 
\end{proof}

%\newpage

\section{Constructing examples}
\label{sec:4}

The aim of this section is to provide 
a concrete combinatorial recipe 
to construct examples of small 
equivariant $G/H$-embeddings with 
Grosshans extensions $H \subset G$.
First, due to our main results,
we obviously have the following general 
recipe.

\begin{construction}
Let $G$ be a connected  
linear algebraic group with 
trivial Picard group 
and only trivial characters.
Every Grosshans extension in $G$ arises 
from the following procedure:
\begin{itemize}
\item
Take a connected Grosshans subgroup 
$F \subset G$ with $\Chi(F) = 0$,
consider the normalizer $N_G(F)$,
and the projection 
$\pi \colon N_G(F) \to N_G(F)/F$.
\item 
Choose a maximal torus $T_F \subset N_G(F)/F$, 
and a surjection $Q \colon \Chi(T_F) \to K$ 
of lattices,
and let $T \subset T_F$ be the corresponding 
subtorus.
\end{itemize}
Then $H_T := \pi^{-1}(T)$
is a Grosshans extension in $G$
with $(H_T)_1 = F$.
The small equivariant $G/H_T$-embeddings
arise from the following procedure:
\begin{itemize}
\item 
Determine the set 
$\Omega(Z)$ of orbit cones of 
the $T$-action on the factorial
affine variety
$Z := \Spec(\mathcal{O}(G)^{F})$.
\item
Fix a 2-maximal collection 
$\Psi \subset \Omega(Z)$
of $T$-orbit cones;
e.g. a collection $\Psi = \Psi_\kappa$ arising 
from a GIT-cone $\kappa \subset K_\QQ$.
\end{itemize} 
Then the 2-maximal collection $\Psi$ defines an 
open subset $U(\Psi) \subset Z$, 
and the quotient $U(\Psi) \quot T$ is 
a small equivariant $G/H_T$-embedding.
\end{construction}

So, the starting point of this construction 
is the choice of a connected Grosshans 
subgroup $F \subset G$;
recall that Remark~\ref{Grosshansexamples}
gives a list of examples.
The second part surely requires a 
certain knowledge of the algebra
$\mathcal{O}(G)^F$; we refer 
to~\cite[Theorem~2]{AT} for 
a detailed study in the case of $G$
being semisimple and $F$ being the 
unipotent 
radical of a parabolic subgroup
$P \subset G$.
 
We now introduce a class of 
Grosshans subgroups
$F \subset G$, the extensions $H_T \subset G$ 
of which allow
a purely combinatorial construction of small
equivariant $G/H_T$-embeddings.

\begin{definition}
\label{suitabledef}
We call a connected Grosshans subgroup
$F \subset G$ {\em suitable} if 
there is a system 
$\{f_1, \ldots, f_r\} \subset \mathcal{O}(G)^F$ 
of $T_F$-homogeneous prime generators
such that every cone in $\Chi_{\QQ}(T_F)$
spanned by some of the weights of 
the $f_i$ is an orbit cone of the 
$T_F$-action on $Z$.
\end{definition}

Here comes the concrete recipe for the
construction of small equivariant embeddings
when starting with a suitable Grosshans 
extension:

\begin{construction}
\label{suitableconstr}
Let $F \subset G$ be a suitable Grosshans 
subgroup, fix a maximal torus 
$T_F \subset N_G(F)/F$, and a system 
$\{f_1, \ldots, f_r\} \subset \mathcal{O}(G)^F$ 
of $T_F$-homogeneous prime generators
as in Definition~\ref{suitabledef}.
\begin{itemize}
\item 
Choose 
a surjection $Q \colon \Chi(T_F) \to K$ 
of lattices,
and let $T \subset T_F$ be the corresponding 
subtorus.
\item
Determine the images $u_i := Q(\deg(f_i))$.
Then the set $\Omega(Z)$ of $T$-orbit cones 
consists of all
$\cone(u_{i_1}, \ldots, u_{i_p})$,
where
$\{i_1, \ldots, i_p\} \subset \{1, \ldots, r\}$.
\item
Fix a 2-maximal collection 
$\Psi \subset \Omega(Z)$
of $T$-orbit cones;
e.g. a collection $\Psi = \Psi_\kappa$ arising 
from a GIT-cone $\kappa \subset K_\QQ$.
%and consider the corresponding 
%2-maximal set $U(\Psi)$.  
\end{itemize} 
Then the 2-maximal collection $\Psi$ defines an 
open subset $U(\Psi) \subset Z$, 
and the quotient $U(\Psi) \quot T$ is 
a small equivariant $G/H_T$-embedding.
\end{construction}

In order to show that this construction
really leads to concrete examples, 
we now present some classes of suitable 
Grosshans subgroups.
We begin with an example
providing spherical varieties.

\begin{proposition}
\label{spherical}
Let $G$ be semisimple simply connected and 
$F \subset G$ be a maximal unipotent subgroup. 
Then $F$ is suitable in $G$.
\end{proposition}

\begin{proof}
Let $\gamma_1,\dots,\gamma_s$ be
fundamental weights of $G$ with respect to a 
Borel subgroup $B=T_FF$, 
and $V(\gamma_1),\dots,V(\gamma_s)$ be corresponding
simple $G$-modules with highest vectors 
$v_{\gamma_i}\in V(\gamma_i)$. 
Then~\cite[Theorem~5.4]{Gr} tells us that 
$$
Z
\; = \; 
\b{G(v_{\gamma_1}, \ldots, v_{\gamma_s})}
\; \subset \; 
V(\gamma_1)\oplus\dots\oplus V(\gamma_s),
$$
and the maximal torus $T_F \cong N_G(F)/F$ 
acts on the variety $Z = \Spec(\mathcal{O}(G)^F)$ 
by means of
$$
t \mal (v_1,\dots,v_s) 
\; = \; 
(\gamma_1(t^{-1})v_1,\dots,\gamma_s(t^{-1})v_s).
$$
For any subset $J\subseteq\{1,\dots,s\}$ consider the
$G$-orbit through a generic point 
$$ 
v_J 
\; \in  \; 
\bigoplus_{j \in J} V(\gamma_j)^F
\; \subset \;
V(\gamma_1)\oplus\dots\oplus V(\gamma_s)
$$ 
The orbit $G \mal v_J$ is contained in $Z$ and
this implies that $F$ is suitable in $G$.
\end{proof}

The following class of examples 
may produce homogeneous spaces of 
arbitrary high complexity.
In particular, they cannot be treated by 
spherical methods, and hence, we go 
a little bit more into detail as before.

\begin{proposition}
\label{xe1}
Let $G := \SL(m)$ act diagonally 
on $(\KK^m)^s$, where $s\le m-1$, 
and consider the isotropy subgroup
$$ 
F 
\; := \;
G_{(e_1, \ldots, e_s)}
\; = \; 
\left\{ 
{\tiny
\left[
\begin{array}{cc}
E_s & A \\
0 & B 
\end{array}
\right]
}; \; 
B \in \SL(m-s), \; A \in \Mat(s \times (m-s)) 
\right\}.
$$
Then $F$ is a connected Grosshans subgroup of $G$,
and a possible maximal torus $T_F \subset N_G(F)/F$ 
is the isomorphic image of  
$$
T'_F 
\; := \; 
\{\text{diag}(t_1,\dots,t_s,t^{-1},1,\dots,1); \; 
t_i \in \KK^*, \; t=t_1 \dots t_s \}
\; \subset \; 
N_G(F). $$
Moreover, we have $Z = \Spec(\mathcal{O}(G)^F) = (\KK^m)^s$,
and the torus $T_F$ acts on the variety $Z$ via 
$$
t \mal (v_1,\dots,v_s) \; = \; (t_1^{-1}v_1,\dots,t_s^{-1}v_s).
$$
In particular, every cone generated by weights of the 
coordinate functions is a $T_F$-orbit cone, and thus
$F$ is suitable in $G$.
\end{proposition}

\begin{proof}
The complement of the open $G$-orbit in $(\KK^m)^s$ 
is the variety of collections of 
linearly dependent vectors, thus it has codimension $\ge 2$. 
This implies
$$
Z \; = \; \Spec(\Of(G/F)) \cong \; (\KK^m)^s.$$ 
In particular, $F$ is a Grosshans subgroup of $G$.
The normalizer $N_G(F)$ coincides 
with the maximal parabolic subgroup
$$
P
\; = \;
\left\{ 
{\tiny
\left[
\begin{array}{cc}
C & A \\
0 & B 
\end{array}
\right]} 
\right\}
\; \subset \;
\SL(m),
$$
and we have $N_G(F)/F\cong\GL(s)$. 
Clearly, the projection $\pi$ maps 
$T'_F$ isomorphically onto
a maximal torus of $\GL(s)$. 
The further statements are obvious. 
\end{proof}

In the setting of Proposition~\ref{xe1},
Construction~\ref{suitableconstr} 
produces small equivariant $G/H$-embeddings 
that come with the structure of a toric 
variety: the action of the torus  
$\TT^{ms}$ on $(\KK^m)^s$ commutes with that 
of $T$, and hence descends to the varieties
$U(\Psi) \quot T$, where $\Psi$ a 2-maximal 
collection, making them into toric varieties.

So it is natural to ask which
further small equivariant 
$G/H$-embeddings have additionaly
a structure of a toric variety.
First of all, 
using~\cite[Cor.~4.5]{BeHa1},
Corollaries~\ref{A2maxsmallemb}
and~\ref{2completechar}
and the fact that
any $T$-maximal subset
of a linear torus actions on 
$\KK^n$ is automatically a
$(T,2)$-set, we obtain
the following.

\begin{remark}
For a toric variety $X$ with free 
divisor class group and 
$\mathcal{O}^*(X) = \KK^*$,
the following statements are equivalent.
\begin{enumerate}
\item 
The variety $X$ is 2-complete.
\item 
The variety $X$ is $A_2$-maximal.
\item 
The fan of $X$ cannot be enlarged 
without adding new rays.
\end{enumerate}
\end{remark}

We show now that besides 
the examples produced via 
Proposition~\ref{xe1} and 
Construction~\ref{suitableconstr},
there is only a very limited list 
of 2-complete small equivariant 
$G/H$-embeddings $X$, where $G$ is
simple and $X$ is toric.

\begin{proposition}
\label{toricclass}
Let $G$ be a simple simply connected
linear algebraic group.
Then the 2-complete small
equivariant $G/H$-embeddings~$X$,
where $H \subset G$ is 
a connected subgroup,
and $X$ admits the structure of
a toric variety,
arise via
Construction~\ref{suitableconstr} 
from the following 
list:
\begin{enumerate}
\item
$G = \SL(m)$ and $Z = (\KK^m)^s$ 
with the diagonal $G$-action and
$F \subset G$ etc. 
as in Proposition~\ref{xe1},
or the dual $G$-module $Z^*$ 
with  the analogous data 
$F \subset G$ etc., 
\item 
the group $G$ and the $G$-modules
$Z$ and their duals $Z^*$, 
where $G$ and $Z$ are
as listed below
$$
\begin{array}{lcl}
G = \SL(2m+1), 
& \quad &
\textstyle
Z = \bigwedge^2\KK^{2m+1},
\\
G = \SL(2m+1), 
& \quad & 
\textstyle
Z = \bigwedge^2\KK^{2m+1} \times \bigwedge^2\KK^{2m+1},
\\
G = \SL(2m+1),
& \quad &   
\textstyle
Z = \bigwedge^2\KK^{2m+1} \times (\KK^{2m+1})^*,
\\
G = \Sp(2m),  
& \quad &   
Z = \KK^{2m},
\\
G = \Spin(10), 
& \quad &
Z = \KK^{16},
\end{array}
$$
and $F = H_1$ is the stabilizer of
a generic point in $Z$ ($Z^*$), 
and the system of generators 
$\{f_1, \ldots, f_r\}$ may be taken as the set 
of coordinate functions.
\end{enumerate}
\end{proposition}

\begin{proof}
Note that $\Cl(X) \cong \Cl(G/H)$
is free, and by~\cite{Co}, the 
toric variety $X$ 
has a polynomial ring as total coordinate ring.
Hence $X$ may be obtained as a good 
$T$-quotient of an open subset 
$U \subset \KK^l$. 
By our assumption, $U$ is a $(T,2)$-maximal
subset of $\KK^l$ and thus, 
by Theorem~\ref{embclass}, 
corresponds to an interior 2-maximal 
collection $\Psi$ of $T$-orbit cones.

Moreover, we infer from 
Proposition~\ref{grossextchar}
and Lemma~\ref{setup} that $H$ is a Grosshans
extension, and that the $G$-action on $X$ 
lifts to a prehomogeneous $G$-action 
on $\KK^l$ commuting
with the $T$-action. 
The $(G\times T)$-action on $\KK^l$ 
is linearizable, see~\cite[Prop.~5.1]{KP} 
for details.
Hence $V$ is a direct sum of simple 
$(G\times T)$-modules, and the $T$-action 
on $V$ is given as in Construction~\ref{suitableconstr}.

Thus, to conclude the proof, we have to 
say what are the possible 
prehomogeneous $G$-modules.
The list of them was obtained in~\cite{Vi}.
\end{proof}

Let us take a quick look at a concrete example.
Assume that $H\subset G$ is a Grosshans extension
such that $G/H$ is quasi-affine. 
If $\Chi(H)=0$, then Theorem~\ref{embclass}
implies that $Z = \Spec(\mathcal{O}(G/H))$
is the only $A_2$-maximal small 
$G/H$-embedding. 
If $H$ has non-trivial characters,
this need no longer hold, as we shall 
see now.

\begin{example}
\label{hyperbolic}
In the notation of Proposition~\ref{xe1}, 
take $m=3$ and $s=2$.
So, we have $G = \SL(3)$ acting diagonally
on $(\KK^3)^2$.
Consider the subtorus
$$
T' 
\; := \; 
\text{diag}(t,t^{-1},1)
\; \subset \; 
T'_F
\; \subset \; 
N_G(F),
$$ 
and set $T = \pi(T')$,
where, as before, 
$\pi \colon N_G(F) \to N_G(F)/F$
is the projection. 
Then the corresponding 
Grosshans extension is 
$$
H_T
\; = \; 
\left\{ 
{\tiny
\left[
\begin{array}{ccc}
t & 0 & a \\
0 & t^{-1} & b \\
0 & 0 & 1 
\end{array}
\right]}; \  t\in\KK^*, \ a,b\in\KK
\right\},
$$
The algebra $\mathcal{O}(G)^F = \mathcal{O}(Z)$ 
is generated by the coordinate functions of 
$Z = (\KK^3)^2$,
and hence the weights of the generators 
in $\ZZ = \Chi(T)$ are $u_1 = 1$ and $u_2 = -1$.
The collection of orbit cones and the possible 
2-maximal collections are given by 
$$
\Omega(Z)
\; = \; 
\{\QQ,\QQ_{\ge 0},\QQ_{\le  0},0\},
\quad 
\Psi_0 
\; = \; 
\{\QQ,0\}, 
\quad 
\Psi_1
\; = \; 
\{\QQ,\QQ_{\ge 0}\},
\quad 
 \Psi_2 = \{\QQ,\QQ_{\le 0}\}. 
$$

Thus, we see that for the homogeneous
space $G/H_T$ there are, up to isomorphism,
precisely three $A_2$-maximal small 
equivariant $G/H_T$-embeddings.
We will discuss them below a little 
more in detail.

The set $U(\Psi_0)$ associated to 
$\Psi_0$ is the whole $Z=(\KK^3)^2$. 
The resulting small equivariant $G/H_T$-embedding
$X_0 = U(\Psi_0) \quot T$ is an affine cone 
with apex $x_1 \in X_0$; it may be 
realized in the 
$G$-module $\KK^3\otimes\KK^3$ 
as the closure 
of the $G$-orbit through $(e_1\otimes e_2)$
with the quotient map 
$$
U(\Psi_0) \; \to \; X_0, \qquad
(v_1,v_2) \; \mapsto \; v_1\otimes v_2.
$$ 

For the collection $\Psi_1$, one has 
$U(\Psi_1)=\{(v_1,v_2); \;  \ v_1\ne 0\}$.
The resulting small equivariant $G/H_T$-emvedding
$X_1 = U(\Psi_1)\quot T$ is quasi-projective
but not affine.
Indeed, the quotient map may be realized via 
$$
U(\Psi_1) \; \to \; (\KK^3\otimes \KK^3) \times \PP^2, 
\qquad  
(v_1,v_2) \; \mapsto \; 
(v_1\otimes v_2, \langle v_1\rangle).
$$
From Theorem~\ref{embclass} we know that
there is a morphism $X_1 \to X_0$ of 
equivariant $G/H_T$-embeddings.
In fact, this is the projection
to $\KK^3\otimes \KK^3$;
this map is an isomorphism
over $X_0 \setminus \{x_1\}$, 
and the fibre over the apex $x_1$
is isomorphic to $\PP^2$. 

The variety $X_2 = U(\Psi_2) \quot T$ 
is isomorphic to $X_1$ as a 
$G$-variety, but not 
as a $G/H_T$-embedding 
(there is no base point preserving 
equivariant morphism).
We may realize $X_2$ by the same construction
as $X_1$ but twisted by the automorphism 
$\theta$ of $Z=\KK^3\oplus\KK^3$,
given by $\theta(v_1,v_2)=(v_2,v_1)$.
\end{example}

In order to see that our construction
also may produce non-toric examples,
look at the following case.

\begin{proposition}\label{xe2}
Let $\KK^{2m}$ be the symplectic vector 
space with the skew-symmetric
bilinear form $\bangle{.\,,.}$, given as
$$ 
{\tiny
\left[
\begin{array}{cc}
0 & E_m \\
-E_m & 0
\end{array}
\right]},
$$
and $G=\Sp(2m)$ be the symplectic group.
Consider the diagonal $G$-action 
on $(\KK^{2m})^s$, where $s \le m$,
and the isotropy group
$$ 
F 
\; := \;
G_{(e_1, \ldots, e_s)}.
$$
Then $F$ is a connected Grosshans subgroup of $G$,
and a possible maximal torus $T_F \subset N_G(F)/F$ 
is the isomorphic image of  
$$
T'_F 
\; := \; 
\{\text{diag}(t_1,\dots,t_s,1,\dots,1,t_1^{-1},\dots,t_s^{-1},1,\dots,1,); 
\; t_i \in \KK^*\}
\; \subset \; 
N_G(F). 
$$
The affine variety $Z = \Spec(\mathcal{O}(G)^F)$ 
can be realized as the $G$-orbit closure of 
$(e_1, \ldots, e_s)$ and is given by 
$$
Z
\; = \; 
\{(v_1,\dots,v_s); \;   \bangle{v_i,v_j} = 0 \ \forall \ i,j \}.
$$
The action of $T_F$ on the variety $Z$ is given as 
$$
t \mal (v_1,\dots,v_s)
\; = \; 
(t_1^{-1}v_1,\ldots,t_s^{-1}v_s).
$$
Every cone generated by weights of the 
restricted coordinate functions is a 
$T_F$-orbit cone, and thus
$F$ is suitable in $G$. 
\end{proposition}

\begin{proof}
First, note that we have $Z = G \mal (L)^s$,
where $L = \langle e_1,\dots, e_m\rangle$ 
is a Lagrangian subspace. 
This shows that the complement of 
the open orbit $G \mal (e_1, \ldots, e_s)$ 
has codimension at least two in $Z$.
Moreover, Serre's Criterion of normality 
shows that $Z$ is normal.
This implies $\mathcal{O}(Z) = \mathcal{O}(G/F)$.

Secondly, the normalizer $N_G(F)$ 
is again a maximal 
parabolic subgroup of $G$, 
we have $N_G(F)/F\cong\GL(s)$, 
and the claim follows.
\end{proof}

\begin{remark}
If we take $s=m$ in the setting of 
Proposition~\ref{xe2}, 
then the subgroup $F$ is the unipotent radical
of the maximal parabolic subgroup $P\subset G$ 
corresponding to the long simple root, 
and is given by 
$$
F 
\; = \;
\left\{
{\tiny
\left[
\begin{array}{cc}
E_m & A \\
0 & E_m
\end{array}
\right]
};  \; 
A=A^T
\right\}.
$$
\end{remark}

\begin{remark}
In Propositions~\ref{spherical}, \ref{xe1} and~\ref{xe2} 
the weights $u_1,\ldots,u_r$ of the 
coordinate functions generate
a regular cone in $\Chi_{\QQ}(T_F)$, 
and the $T_F$-orbit cones are precisely 
the faces of this cone. 
\end{remark}

To finish the discussion on suitable subgroups, 
we give an example showing that not any connected Grosshans
subgroup $F\subset G$ with $\Chi(F)=0$ is suitable.

\begin{example}
Let $F$ be a connected semisimple subgroup of $G$. 
Then $Z=G/F$ and the only orbit cone for the $T_F$-action 
on $Z$ is $\Chi_{\QQ}(T_F)$. 
This shows that $F$ is not suitable in $G$.
\end{example}

%\newpage

\section{Geometric properties}
\label{sec:5}

In this section, we show that the 
language of bunched rings developed 
in~\cite{BeHa1}, applies to $A_2$-maximal 
small equivariant $G/H$-embeddings $X$,
provided that $H \subset G$ is a Grosshans
extension.
This enables us to study basic geometric
properties of $X$.
For example, we obtain existence of 
projective small equivariant $G/H$-embeddings
with at most $\QQ$-factorial singularities,
and we can easily produce 
homogeneous spaces $G/H$ that do not
admit any smooth small equivariant  
completion.

Let us briefly recall the concepts of~\cite{BeHa1}.
In the sequel, $R$ denotes a factorial, 
finitely generated $\KK$-algebra, 
faithfully graded by some lattice $K \cong \ZZ^{k}$
such that $R^* = \KK^*$ holds.
Here, faithfully graded means that $K$ is 
generated as a lattice by the degrees $w \in K$ admitting 
nontrivial homogeneous elements $f \in R_w$.

Moreover, 
$\mathfrak{F} =\{f_{1}, \dots, f_{r}\} \subset R$ is a 
system of homogeneous pairwise non associated 
nonzero prime elements generating $R$ as an algebra.
Note that due to $R^* = \KK^*$ such systems always
exist.
Since we assume the grading to be faithful,
the degrees $\deg(f_{i})$ generate the
lattice $K$.

The {\em projected cone $(E \topto{Q} K, \gamma)$
associated\/} to the system of generators 
$\mathfrak{F} \subset R$ consists of the surjection 
$Q$ of the lattices $E := \ZZ^{r}$ and $K$ sending 
the $i$-th canonical base vector $e_{i} \in \ZZ$ 
to the degree $\deg(f_{i}) \in K$, 
and the cone $\gamma \subset  E_{\QQ}$ generated by 
$e_{1}, \dots, e_{r}$.

\begin{definition}
\label{fdefs}
Let $(E \topto{Q} K, \gamma)$ be
the projected cone associated to
$\mathfrak{F} \subset R$, and suppose that
for each facet $\gamma_{0} \preceq \gamma$,
the image $Q(\gamma_{0} \cap E)$ generates 
the lattice $K$.
\begin{enumerate}    
\item \label{fface}    
A face $\gamma_{0} \preceq \gamma$ is called an
{\em $\mathfrak{F}$-face\/} if the product
over all $f_{i}$ with $e_{i} \in \gamma_{0}$
does not belong to the ideal
$\sqrt{\bangle{f_{j}; \; e_{j} \not\in \gamma_{0}}} \subset R$.
\item \label{fbunch}
An {\em $\mathfrak{F}$-bunch\/} is  a nonempty 
collection $\Phi$ of projected $\mathfrak{F}$-faces 
with the following properties:
\begin{itemize}
\item 
a projected $\mathfrak{F}$-face $\tau$ belongs to $\Phi$ if 
and only if for each $\tau \neq \sigma \in \Phi$ we have
$\emptyset \neq \tau^{\circ} \cap \sigma^{\circ} \neq 
\sigma^{\circ}$,
\item 
for each facet $\gamma_{0} \preceq \gamma$,
there is a cone $\tau \in \Phi$ such that 
$Q(\gamma_{0})^{\circ} \supset \tau^{\circ}$ holds.
\end{itemize}
\end{enumerate}
If $\Phi$ is an $\mathfrak{F}$-bunch
in the projected cone $(E \topto{Q} K, \gamma)$ 
associated to $\mathfrak{F} \subset R$, 
then the triple $(R,\mathfrak{F},\Phi)$
is called a {\em bunched ring}. 
\end{definition}

Now consider the affine variety $Z := \Spec(R)$,
the torus $T := \Spec(\KK[K])$, and the action
$T \times Z \to Z$ given by the $K$-grading of
$R$.
The following statements put the above definitions
into a more geometric framework. 

\begin{lemma}
\label{fdefsgeo}
Let $(E \topto{Q} K, \gamma)$ 
be the projected cone 
associated to $\mathfrak{F} \subset R$.
Then the following statements hold:
\begin{enumerate}
\item 
The projected $\mathfrak{F}$-faces
are precisely the orbit cones of the 
$T$-action on $Z$.
\item 
There is a canonical injection
\begin{eqnarray*}
\{\mathfrak{F}  \text{-bunches} \}
& \rightarrow &
\{\text{2-maximal collections in } \Omega(Z)\}
\\
\Phi
& \mapsto &
\left\{
\omega(z); \; 
z \in Z, \; 
\tau^\circ \subset \omega(z)^\circ 
\text{ for some } \tau \in \Phi
\right\}.
\end{eqnarray*}
\end{enumerate}
\end{lemma}

\begin{proof}
To prove the first statement,
we first note that the defining condition
of an $\mathfrak{F}$-face has the following
geometric meaning: it says that 
$\gamma_0 \preceq \gamma$ is an 
$\mathfrak{F}$-face if and only if 
there is a point $z \in Z$ such that 
\begin{eqnarray}
\label{Fgeom}
f_i(z) \; \ne  \; 0
& \Leftrightarrow &
e_i \in \gamma_0.
\end{eqnarray}

Consider any orbit cone $\omega(z)$
where $z \in Z$. We claim that 
$\omega(z) = Q(\gamma_0)$ holds
for the face $\gamma_0 \preceq \gamma$
defined by
\begin{eqnarray}
\label{orbitcone2fi}
\gamma_0
& = & 
\cone(e_i; \; f_i(z) \ne 0)
\end{eqnarray}
Obviously, we have $Q(\gamma_0) \subset \omega(z)$.
For the converse inclusion, consider any 
homogeneous $h \in R$ 
with $h(z) \ne 0$.
Then we have a representation
\begin{eqnarray*}
h
& = & 
\sum \alpha_\nu f_1^{\nu_1} \ldots f_r^{\nu_r}
\end{eqnarray*}
with coefficients $\alpha_\nu \in \KK$.
Consequently, the degree of $h$ is a 
positive combination of some of the degrees
of the $f_i$. This shows 
$\omega(z) \subset Q(\gamma_0)$.

Now, given any orbit cone $\omega(z)$,
this cone is 
the image of the face $\gamma_0 \preceq \gamma$
given as in~(\ref{orbitcone2fi}).
Moreover, the point $z$ satisfies~(\ref{Fgeom}), 
showing that $\gamma_0$ is an $\mathfrak{F}$-face.
Conversely, given any $\mathfrak{F}$-face 
$\gamma_0 \preceq \gamma$, consider 
$z \in Z$ as in~(\ref{Fgeom}). 
Then~(\ref{orbitcone2fi}) shows that 
$Q(\gamma_0)$ is the orbit cone of $z$.

The second assertion is simply due to the 
observation that the bunch $\Phi$ may be 
reconstructed from its associated 
2-maximal collection $\Psi$ by
taking the set-theoretically minimal 
elements of $\Psi$.     
\end{proof}

Now suppose, we are in the setting 
of the preceding section. That means 
that $G$ is a connected affine algebraic
group with $\Chi(G) = 0$ and $\Pic(G) = 0$, 
and $H \subset G$, is a Grosshans
extension.
Then $W = G/H_1$ is a quasiaffine 
variety with a finitely generated
algebra $R := \mathcal{O}(W)$ 
of global functions satisfying 
$R^* = \KK^*$.

Moreover, denoting by $K := \Chi(H)$ the
character lattice of $H$, we have the 
canonical action of the torus 
$T = \Spec(\KK[K])$ on the 
quasiaffine variety $W$.
The actions of $G$ and $T$ extend to 
the factorial affine variety 
$Z := \Spec(R)$; in particular the 
$T$-action defines a grading:
\begin{eqnarray*}
R
& := &
\bigoplus_{u \in K} R_{u}.
\end{eqnarray*}
The basic observation of this section
is that, with the above data, we are 
in the setting of Definition~\ref{fdefs}.
More precisely, we observe the following.

\begin{proposition}
\label{smallGH2bunchedring}
There is a system 
$\mathfrak{F} = \{f_1, \ldots, f_r\} \subset R$
of homogeneous pairwise nonassociated 
prime generators with the following properties.
\begin{enumerate}
\item 
Let $(E \topto{Q} K, \gamma)$ be the projected 
cone associated to $\mathfrak{F} \subset R$.
Then for every facet $\gamma_0 \prec \gamma$,
the image $Q(E \cap \gamma_0)$ generates $K$ 
as a lattice.
\item 
We have canonical bijections, inverse to each
other:
\begin{eqnarray*}
\{\mathfrak{F}  \text{-bunches} \}
& \longleftrightarrow &
\{\text{interior 2-maximal collections}\}
\\
\Phi
& \mapsto &
\left\{
\omega(z); \; 
%z \in Z, \; 
\tau^\circ \subset \omega(z)^\circ 
\text{ with } \tau \in \Phi
\right\}
\\
\{\omega \text{ minimal in } \Psi\}
& \mapsfrom & 
\Psi
\end{eqnarray*}
\end{enumerate}
If  $\mathcal{O}(G/H) = \KK$ 
holds, then even every system 
$\mathfrak{F} = \{f_1, \ldots, f_r\} \subset R$
of homogeneous pairwise nonassociated 
prime generators fullfills the above
properties.
\end{proposition}

\begin{proof}
In a first step we show that any system
$\mathfrak{F} = \{f_1, \ldots, f_r\} \subset R$
of homogeneous pairwise nonassociated
prime generators satisfies condition~(i).
After suitably renumbering, we may assume 
$\gamma_0 = \cone(e_1, \ldots, e_{r-1})$.
Consider any $z \in Z$ with $f_r(z) = 0$.
Then the isotropy group $T_z \subset T$ has 
as its character group $K/K(z)$, where  
\begin{eqnarray*}
K(z)
& := & 
\bangle{u \in K; \; f(z) \ne 0 \text{ for some } f \in R_u} 
\\
& \subset &
\bangle{\deg(f_1), \ldots, \deg(f_{r-1})}.
\end{eqnarray*}
Since the zero set $V(Z,f_r)$
must intersect the big $G$-orbit $W = G/H_1$,
and $T$ acts freely on this set, we see
there are points $z \in V(Z,f_r)$ with $K(z) = K$.
Thus, the displayed formula shows that 
the image $Q(\gamma_0 \cap E)$
generates $K$ as a lattice.

If we are in the general case, i.e., 
$\mathcal{O}(G/H)$ may contain nonconstant
functions, then we have to provide a suitable 
system $\mathfrak{F} \subset R$ of generators.
For this, we first take any collection of 
homogeneous 
pairwise nonassociated prime elements 
$\{f_1, \ldots, f_r\}$ generating $R$.
For each of these $f_i$, we consider the 
$G$-stable vector subspace $V_i \subset R$ 
generated by $G \mal f_i$.

Since $G$ acts with an open orbit on $Z$,
we have $R^G = \KK$, and hence $V_i$ is 
a nontrivial $G$-module.
Since we have $\Chi(G) = 0$, we even see 
that $\dim(V_i) > 1$ holds.
So, there is an $f'_i = g_i \mal f_i$,
which is not proportional to $f_i$.
Since we have $R^* = \KK^*$ this even means
that $f_i$ and $f_i'$ are nonassociated
primes.

Adding appropriate elements $f_i'$, 
we may enlarge the initial system
of generators
$\{f_1, \ldots, f_r\}$ such that for any
$i$, there is a $j \ne i$ with 
$\deg(f_i) = \deg(f_j)$.
Then it is obvious that this new 
complemented system satisfies
\begin{equation}
\label{inner}
Q(\gamma_0) \; = \; Q(\gamma) \; = \; \omega(Z)
\text{ for every facet } \gamma_0 \preceq \gamma.
\end{equation}

Let us show that this property
gives the second condition.
By Lemma~\ref{fdefsgeo}, we have a
canonical injection from the 
$\mathfrak{F}$-bunches 
to the 2-maximal collections.
Our task is to show that the
image consists precisely of 
the interior 2-maximal collections.

Given an $\mathfrak{F}$-bunch $\Phi$,
condition~(\ref{inner}) implies
that $\omega(Z)$ occurs in the 
associated 2-maximal collection $\Psi$.
Since $\omega(Z)$ is the generic orbit 
cone, $\Psi$ is an interior collection.
Conversely, given an interior 2-maximal 
collection, (\ref{inner}) yields 
that the corresponding collection $\Phi$ 
of its set-theoretically minimal cones 
satisfies~\ref{fdefs}~(ii), and hence is an 
$\mathfrak{F}$-bunch.

Now suppose that $\mathcal{O}(G/H) = \KK$ holds.
Then we have to show that any system
$\mathfrak{F} = \{f_1, \ldots, f_r\} \subset R$
of homogeneous pairwise nonassociated
prime generators satisfies condition~(\ref{inner}).
First note that $\mathcal{O}(G/H) = \KK$ 
implies that $Q(\gamma) = \omega(Z)$ is pointed.
Thus, in order to obtain~(\ref{inner}),
we have to show that each extremal
ray of $\omega(Z)$ contains at least
two of the $u_i = \deg(f_i)$.

Let us verify this.  Clearly, 
$\varrho$ contains at least one 
$u_i$.
In order to see that there must 
be a second one, consider any nontrivial 
translate $h := g \mal f_i$.
This is as well a homogeneous function 
of degree $u_i$. 
Since we have $\Chi(G) = 0$ the elements 
$f_i$ and $h$ are linearly independent. 
%A fortiori, $R^* = \KK^*$, implies that 
%they  are even nonassociated primes.
Thus, there must be a representation 
\begin{eqnarray*}
h
& = & 
\sum \alpha_\nu f_1^{\nu_1} \ldots f_r^{\nu_r}
\end{eqnarray*}
with coefficients $\alpha_\nu \in \KK$ such that
$\alpha_\nu \ne 0$ for at least some 
$\nu$ admitting a $\nu_j \ne 0$ with $j \ne i$.
Since $u_i$ belongs to an extremal ray, 
we obtain $u_j \in \varrho$.
This establishes condition~(\ref{inner}) 
for $\{f_1, \ldots, f_r\}$, from which we 
deduce, as before,
condition~(ii) of the assertion.
\end{proof}

This proposition shows that
for a suitable system 
$\mathfrak{F} = \{f_1, \ldots, f_r\}$
of generators of the algebra 
$R=\mathcal{O}(G/H_1)$, 
the small equivariant $G/H$-embeddings
are precisely the varieties 
\begin{eqnarray*}
X(R,\mathfrak{F},\Phi)
& = & 
U(\Psi) \quot T
\end{eqnarray*}
arising from the bunched rings 
$(R,\mathfrak{F},\Phi)$,
where $\Psi$ is the 2-maximal collection
associated to the $\mathfrak{F}$-bunch $\Phi$.
So, we may apply the results obtained
in~\cite{BeHa1} to describe geometric
properties of $X$.

Let us briefly provide the necessary
notions.
Call an $\mathfrak{F}$-face 
$\gamma_{0} \preceq \gamma$ {\em relevant\/}
if $Q(\gamma_{0})^{\circ} \supset \tau^{\circ}$ 
holds for some $\tau \in \Phi$,
and denote by $\rlv(\Phi)$ the collection
of relevant $\mathfrak{F}$-faces.
The {\em covering collection} of $\Phi$ 
is the collection $\cov(\Phi) \subset \rlv(\Phi)$
of set-theoretically minimal members of 
$\rlv(\Phi)$. 
Here are some of the results of~\cite{BeHa1}.

\begin{theorem}
\label{bunchedrings}
For a suitable choice of 
$\mathfrak{F} \subset R$,
let $X := X(R,\mathfrak{F},\Phi)$ 
be the small $G/H$-embedding 
arising from an $\mathfrak{F}$-bunch 
$\Phi$. 
Then the following statements 
hold:
\begin{enumerate}
\item 
The variety 
$X$ is locally factorial if and only if
$Q(\gamma_0 \cap E)$ generates the lattice 
$K$ for every $\gamma_0 \in \rlv(\Phi)$.
\item 
The variety 
$X$ is $\QQ$-factorial if and only if 
any cone of $\Phi$ is of full dimension 
in $K_\QQ$. 
\item 
The rational divisor class group of $X$ is given by
$\Cl_\QQ(X) \cong K_\QQ$, and inside $K_\QQ$
the Picard group of $X$ is given by 
\begin{eqnarray*}
\Pic(X) 
& = & 
\bigcap_{\gamma_{0} \in \cov(\Phi)} Q(\lin(\gamma_{0}) \cap E).
\end{eqnarray*}
%\item 
Moreover, inside $K_\QQ$, the cones of semiample and ample 
divisor classes of X are given by
$$
\begin{array}{ccc}
\displaystyle \SAmple(X) \; = \; \bigcap_{\tau \in \Phi} \tau,
& \qquad &
\displaystyle \Ample(X) \; = \; \bigcap_{\tau \in \Phi} \tau^{\circ}.
\end{array}
$$
\end{enumerate}
\end{theorem}

Note that the semiample and the ample 
cone depend only on the bunch $\Phi$,
and might as well be expressed in terms
of the corresponding 2-maximal 
collection~$\Psi$.

As a very first application of 
Theorem~\ref{bunchedrings}, we give
an existence statement on small
equivariant $G/H$-embeddings
in the spirit of~\cite[Th\'{e}or\`{e}me~1]{BBII},
but additionally guaranteing
``mild'' singularities.

\begin{corollary}
\label{existence}
Let $G$ be a connected 
linear algebraic group 
trivial Picard group and
only trivial characters, 
and let $H \subset G$
be a Grosshans extension with 
$\mathcal{O}(G)^H = \KK$.
Then there exists a projective
small equivariant $G/H$-embedding
with at most $\QQ$-factorial singularities. 
\end{corollary}

\begin{proof}
Let $R := \mathcal{O}(G)^{H_1}$,
let $\mathfrak{F} \subset R$ be any 
system of pairwise nonassociated homogeneous
prime generators, and take any $\mathfrak{F}$-bunch 
arising from a GIT-cone of full dimension.
Then the corresponding small $G/H$-embedding 
is as wanted.
\end{proof}

We will now apply the language of bunched rings 
to study the concrete examples arising from 
the constructions of the preceding section more 
in detail.
As it concerns a good part of them, we first note
the following.

\begin{remark}
\label{toric}
Given a bunched ring $(R,\mathfrak{F},\Phi)$,
where $R = \KK[z_1, \ldots, z_r]$ is a polynomial
ring, and $\mathfrak{F} = \{z_1, \ldots, z_r\}$
consists of the indeterminates, we are in the
setting of toric varieties, and then, in addition
to Theorem~\ref{bunchedrings}, there are 
simple combinatorial 
criteria for smoothness~\cite[Prop.~8.3]{BeHa0} 
and completeness~\cite[Prop.~8.6]{BeHa0}.
\end{remark}

Now we begin the study of examples.
The first one gives smooth small 
equivariant $G/H$-embeddings;
recall from~\cite[Sec.~5]{BBK} 
that the existence of a smooth projective
small $G/H$-embedding implies that 
$G/H$ is generically rationally connected.

\begin{example}
\label{smoothemb}
In the setting of 
Proposition~\ref{xe1}, 
let $m:=4$ and $s:=3$.
Then the torus 
$T_F \subset N_G(F)/F$
is of dimension three, 
and it may be identified
with
$$
\{\text{diag}(t_1,t_2,t_3,t_1^{-1}t_2^{-1}t_3^{-1}); \; 
t_i \in \KK^* \}
\; \subset \; 
N_G(F). $$ 
We consider the projection
$\Chi(T_F) \to \ZZ^2$
sending the canonical generators
of the character group $\Chi(T_F)$ 
to the following 
lattice vectors
$$ 
u_1 \; := \; (1,0),
\qquad
u_0 \; := \; (1,1),
\qquad
u_2 \; := \; (0,1).
$$
Thus, speaking more concretely, we deal with 
$G = \SL(4)$, the Grosshans subgroup 
$F \subset G$ as in~\ref{xe1}, 
a twodimensional torus
$T \subset N_G(F)/F$, and the Grosshans 
extension $H_T \subset G$ given by
$$
H_T
\; = \; 
\left\{ 
{\tiny
\left[
\begin{array}{cccc}
t_1 & 0 & 0 & a_1     \\
0 & t_1t_2 & 0 & a_2 \\
0 & 0 & t_2  & a_3 \\
0 & 0 & 0 & t_1^{-2}t_2^{-2} 
\end{array}
\right]}; \;  t_1, t_2 \in\KK^*, \; a_i \in\KK
\right\},
$$

Let us determine the bunched rings 
$(R,\mathfrak{F},\Phi)$ describing 
the possible $A_2$-maximal small 
equivariant $G/H_T$-embeddings.
First of all, $R = \mathcal{O}(G)^{F}$ 
is the ring of functions of 
$G/F = (\KK^4)^3$. So, $R$ is a polynomial 
ring, and as system of 
generators 
$\{f_1, \ldots, f_{12}\} \subset R$, 
one may take the collection of 
indeterminates.

The remaining task is to determine the possible
$\mathfrak{F}$-bunches.
As we know from Corollary~\ref{Z2proj}, these 
bunches correspond to projective varieties,
and hence we only need to know the GIT-fan
in $\ZZ^2 \cong \Chi(T)$ of the action of $T$ on $G/F$.
This in turn is easy to determine; 
it looks as follows:
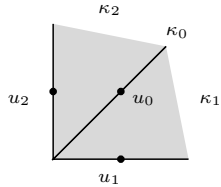
\begin{center}
\begin{picture}(0,0)%
\includegraphics{smoothemb.pstex}%
\end{picture}%
\setlength{\unitlength}{1243sp}%
\begingroup\makeatletter\ifx\SetFigFont\undefined%
\gdef\SetFigFont#1#2#3#4#5{%
  \reset@font\fontsize{#1}{#2pt}%
  \fontfamily{#3}\fontseries{#4}\fontshape{#5}%
  \selectfont}%
\fi\endgroup%
\begin{picture}(3855,3825)(-14,-3334)
\put(2476,-1636){\makebox(0,0)[lb]{\smash{{\SetFigFont{7}{8.4}{\familydefault}{\mddefault}{\updefault}{\color[rgb]{0,0,0}$u_0$}%
}}}}
\put(  1,-1636){\makebox(0,0)[lb]{\smash{{\SetFigFont{7}{8.4}{\familydefault}{\mddefault}{\updefault}{\color[rgb]{0,0,0}$u_2$}%
}}}}
\put(1801,-3211){\makebox(0,0)[lb]{\smash{{\SetFigFont{7}{8.4}{\familydefault}{\mddefault}{\updefault}{\color[rgb]{0,0,0}$u_1$}%
}}}}
\put(3151,-286){\makebox(0,0)[lb]{\smash{{\SetFigFont{7}{8.4}{\familydefault}{\mddefault}{\updefault}{\color[rgb]{0,0,0}$\kappa_0$}%
}}}}
\put(3826,-1636){\makebox(0,0)[lb]{\smash{{\SetFigFont{7}{8.4}{\familydefault}{\mddefault}{\updefault}{\color[rgb]{0,0,0}$\kappa_1$}%
}}}}
\put(1801,164){\makebox(0,0)[lb]{\smash{{\SetFigFont{7}{8.4}{\familydefault}{\mddefault}{\updefault}{\color[rgb]{0,0,0}$\kappa_2$}%
}}}}
\end{picture}%

\end{center}
So, the possible $\mathfrak{F}$-bunches are 
those arising from the interior GIT-cones
$\kappa_1$, $\kappa_0$ and $\kappa_2$ as 
indicated above, and they are explicitly 
given by 
$$ 
\Phi_1 \; = \; \{\kappa_1\},
\qquad  
\Phi_0 \; = \; \{\kappa_0\},
\qquad  
\Phi_2\; = \; \{\kappa_2\}.
$$

Let $X_i$ denote the small equivariant 
$G/H_T$-embedding corresponding to the 
$\mathfrak{F}$-bunch $\Phi_i$.
Then, applying the results on the 
geometry of $X_i$ mentioned in~\ref{bunchedrings}
and~\ref{toric}, 
we see for example that $X_1$ and $X_2$ 
are smooth, whereas $X_0$ has 
non-$\QQ$-factorial singularities.

Moreover, all varieties $X_i$ have a divisor class group of
rank two, and $X_0$ has a Picard group 
of rank one.   
Finally, Theorem~\ref{embclass} tells
us that the possible morphisms
of $G/H_T$-embeddings are
$$ 
X_1 \; \longrightarrow \; X_0 \; \longleftarrow \; X_2.
$$  
By determining explicitly the varieties
$U(\kappa_i)$ over $X_i$ one may describe 
these morphisms explicitly, and it turns
out that for $i=1,2$ the exceptional locus 
of $X_i \to X_0$ is isomorphic to 
$\PP^3 \times \PP^3$ and is contracted
to a $\PP^3$ lying in $X_0$. 

Moreover, one obtains that, as $G$-varieties, 
$X_1$ and $X_2$ are
isomorphic, but, of course, not as $G/H$-embeddings.
This shows in particular that $N_G(H)/H$ is not
contained in the group of $G$-equivariant 
automorphisms of $X_1$.
\end{example}

By slight modification of the preceding
example, we present a homogeneous space 
$\SL(4)/H$ admitting equivariant
completions with small boundary
but no smooth ones.
Existence of such examples is due to 
M.~Brion, as mentioned in~\cite{BBK}.

\begin{example}
\label{nosmoothemb}
In the setting of 
Proposition~\ref{xe1}, 
let $m:=4$ and $s:=3$.
As before, $T_F$ is of dimension three, 
but now we consider the projection
$\Chi(T_F) \to \ZZ^2$
sending the canonical generators
to
$$ 
u_1 \; := \; (1,0),
\qquad
u_0 \; := \; (2,3),
\qquad
u_2 \; := \; (0,1).
$$
Concretely this means that we have again 
the Grosshans subgroup $F \subset G = \SL(4)$, 
as in~\ref{xe1}, but another twodimensional
torus $T \subset N_G(F)/F$. The Grosshans 
extension $H_T \subset G$ this time 
is given by
$$
H_T
\; = \; 
\left\{ 
{\tiny
\left[
\begin{array}{cccc}
t_1 & 0 & 0 & a_1     \\
0 & t_1^2t_2^3 & 0 & a_2 \\
0 & 0 & t_2  & a_3 \\
0 & 0 & 0 & t_1^{-3}t_2^{-4} 
\end{array}
\right]}; \;  t_1, t_2 \in\KK^*, \; a_i \in\KK
\right\},
$$
The possible small equivariant 
$G/H$-completions arise from $T$-maximal 
open subsets of $Z = (\KK^4)^3$.
All of them are toric, 
hence $A_2$-maximal,
hence projective, use e.g. Corollary~\ref{Z2proj}.
Thus the GIT-fan for the $T$-action on 
$Z = (\KK^4)^3$ gives the full information;
it look as follows:
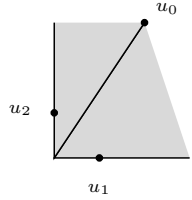
\begin{center}
\begin{picture}(0,0)%
\includegraphics{singemb.pstex}%
\end{picture}%
\setlength{\unitlength}{1243sp}%
\begingroup\makeatletter\ifx\SetFigFont\undefined%
\gdef\SetFigFont#1#2#3#4#5{%
  \reset@font\fontsize{#1}{#2pt}%
  \fontfamily{#3}\fontseries{#4}\fontshape{#5}%
  \selectfont}%
\fi\endgroup%
\begin{picture}(3648,4026)(-14,-3559)
\put(1576,-3436){\makebox(0,0)[lb]{\smash{{\SetFigFont{7}{8.4}{\familydefault}{\mddefault}{\updefault}{\color[rgb]{0,0,0}$u_1$}%
}}}}
\put(  1,-1861){\makebox(0,0)[lb]{\smash{{\SetFigFont{7}{8.4}{\familydefault}{\mddefault}{\updefault}{\color[rgb]{0,0,0}$u_2$}%
}}}}
\put(2926,164){\makebox(0,0)[lb]{\smash{{\SetFigFont{7}{8.4}{\familydefault}{\mddefault}{\updefault}{\color[rgb]{0,0,0}$u_0$}%
}}}}
\end{picture}%

\end{center}
Using Theorem~\ref{bunchedrings}~(i), we see that
each of the three possible projective small 
equivariant $G/H$-embeddings is singular;
in two cases, we have $\QQ$-factorial singularities,
in the remaining one, it is even worse.
\end{example}

\begin{example}
In the setting of Proposition~\ref{xe1},
let $m=7$ and $s=6$.
So, we have $G= \SL(7)$ acting diagonally
on $(\KK^7)^6$, and the torus $T_F$ 
is of dimension six.
Set $K := \ZZ^3$ and consider the map 
$\Chi(T_F) \to K$ sending the 
canonical generators to
$$ 
e_1,
\qquad
e_2,
\qquad
e_3,
\qquad
w_1 \; :=  \; e_1 + e_2,
\qquad
w_2 \; :=  \; e_1 + e_3,
\qquad
w_3 \; :=  \; e_2 + e_3.
$$
Let $\mathfrak{F} = \{f_1, \ldots, f_{42}\}$
be the indeterminates of the polynomial ring
$\mathcal{O}(G)^H$.
Then the following cones 
in $K_\QQ$ define an $\mathfrak{F}$-bunch:
$$ 
\cone(e_3,w_1,w_2),
\qquad
\cone(e_1,w_1,w_3),
\qquad
\cone(e_2,w_2,w_3),
\qquad
\cone(w_1,w_1,w_2),
$$
Combining~\cite[Example~11.2]{BeHa0} and
\cite[Construction~11.4]{BeHa0} shows
that the corresponding small equivariant 
$G/H_T$-embedding $X(R,\mathfrak{F},\Phi)$
is complete and $\QQ$-factorial
but not projective.
\end{example}

Finally, we give a concrete example showing 
that the language of bunched rings also 
applies in the non-toric case. 

\begin{example}
In the setting of Proposition~\ref{xe2}, 
let $m = 3$, and $s=3$.
Then the maximal torus $T_F \subset N_G(F)/F$
is of dimension three.
Consider a surjection $\Chi(T_F) \to \ZZ^2$,
sending the canonical generators to 
$u_1, u_0, u_2 \in \ZZ^2$.
Then our Grosshans extension 
$H_T \subset \Sp(6)$ consists of the 
matrices
$$
%\left\{
{\tiny
\left[
\begin{array}{cccccc}
\chi^{u_1}(t) & 0 & 0 & a_{11} & a_{12} & a_{13} \\
0 & \chi^{u_2}(t) & 0 & a_{12} & a_{22} & a_{23} \\
0 & 0 & \chi^{u_3}(t) & a_{13} & a_{23} & a_{33} \\
0 & 0 & 0      & \chi^{-u_1}(t) & 0 & 0  \\    
0 & 0 & 0      & 0 & \chi^{-u_2}(t) & 0  \\
0 & 0 & 0      & 0 & 0 & \chi^{-u_3}(t)
\end{array}
\right]
}, 
\quad 
%\text{where } 
t \in T; \; 
a_{ij} \in \KK.
%\right\}
$$ 

Taking $u_i$ as in Example~\ref{smoothemb}, 
we obtain three projective small equivariant 
$G/H_T$-embeddings $X_1$, $X_0$, and $X_2$. 
By Theorem~\ref{bunchedrings}, 
the varieties $X_1$ and $X_2$ are locally 
factorial and $X_3$ is not $\QQ$-factorial

In fact, an analysis of the singular locus 
of the cone $Z = \Spec(\mathcal{O}(G)^{F})$ 
shows that open subsets $U(\Psi_1)$ and 
$U(\Psi_2)$ lying over $X_1$ and $X_2$
respectively are smooth. 
Thus~\cite[Proposition~5.6]{BeHa1}
yields that the varieties $X_1$ and $X_2$ 
are even smooth.

Taking $u_i$ as in Example~\ref{nosmoothemb}, 
one obtains another torus $T$, and other 
projective small equivariant 
$G/H_T$-embeddings $X_1$, $X_0$, and $X_2$.
Then Theorem~\ref{bunchedrings} tells us that
$X_1$ and $X_2$ are $\QQ$-factorial, but not
locally factorial.
\end{example}

%\newpage


\begin{thebibliography}{}%
%
\bibitem{AT} I.V.~Arzhantsev, D.A.~Timashev: 
On the canonical embeddings of certain homogeneous spaces.
In ``Lie Groups and Invariant Theory: A.L.~Onishchik's jubilee volume'' 
(E.B.~Vinberg, Editor),
AMS Translations, Series~2, vol.~213, 63--83 (2005)
%
\bibitem{BeHa0} F.~Berchtold, J.~Hausen: 
Bunches of cones in the divisor class group -- 
a new combinatorial language for toric varieties. 
Inter. Math. Research Notices 6, 261--302 (2004)
%
\bibitem{BeHa1} F.~Berchtold, J.~Hausen: Cox rings and combinatorics.
To appear in Transactions of the AMS, {\tt math.AG/0311105}
%
\bibitem{BeHa2} F.~Berchtold, J.~Hausen: GIT-equivalence beyond 
the ample cone,
Preprint, {\tt math.AG/0503107}
%
\bibitem{BB1}
A.~Bia\l ynicki-Birula: Finiteness of the number of maximal 
open subsets with good quotients.  
Transform. Groups  3  (1998),  no. 4, 301--319.
%
\bibitem{BB2}
A.~Bia\l ynicki-Birula: Algebraic quotients. 
In: Encyclopaedia of Mathematical Sciences, 131. 
Invariant Theory and Algebraic Transformation Groups, II. 
Springer-Verlag, Berlin, 2002.
%
\bibitem{BBSw}
A.~Bia\l ynicki-Birula, J.~\'Swi\c{e}cicka:
A recipe for finding open subsets of vector 
spaces with a good quotient.  
Colloq. Math.  77  (1998),  no. 1, 97--114.
%
\bibitem{BBI} F.~Bien, A.~Borel: 
Sous-groupes \'epimorphiques de groupes lin\'eaires alg\'ebriques I.
C. R. Acad. Sci. Paris, t.~315, S\'erie I, 649--653 (1992)
%
\bibitem{BBII} F.~Bien, A.~Borel: 
Sous-groupes \'epimorphiques de groupes lin\'eaires alg\'ebriques II.
C. R. Acad. Sci. Paris, t.~315, S\'erie I, 1341--1346 (1992)
%
\bibitem{BBK} F.~Bien, A.~Borel, J.~Kollar: 
Rationally connected homogeneous spaces.
Invent. Math. 124, 103--127 (1996)
%
\bibitem{BrPr} Brion.~M.; Procesi, C.:
Action d'une tore dans une vari\'{e}t\'{e} projective.
In: Connes, A. et al. (Eds): 
Operator Algebras, Unitary Representations, 
Enveloping Algebras, and Invariant Theory. 
Progress in Mathematics, Vol. 192, 
Birkh\"{a}user, Basel 1990, 509-539
%
\bibitem{Co} D.~Cox: The homogeneous coordinate ring of a toric variety.
J. Alg. Geom. 4, 17--50 (1995)
%
\bibitem{DoHu}
I.V.~Dolgachev, Y.~Hu: Variation of geometric invariant theory quotients. 
With an appendix by Nicolas Ressayre.  
Inst. Hautes Études Sci. Publ. Math.  No. 87 (1998), 5--56.
%
%\bibitem{Ei} M.~Eikelberg: The Picard group of a compact toric variety.
%Results in Math. 22, 509--527 (1992)
%
\bibitem{ElKuWa} J.~Elizondo, K.~Kurano, K.~Watanabe: 
The total coordinate ring of a normal projective variety. 
J. Algebra 276, 625--637 (2004) 
%
\bibitem{Gr} F.D.~Grosshans, Algebraic Homogeneous Spaces 
and Invariant Theory, 
LNM 1673, Springer-Verlag Berlin (1997)
%
\bibitem{Ha1} J.~Hausen: A general Hilbert-Mumford criterion.  
Ann. Inst. Fourier (Grenoble)  53  (2003),  no. 3, 701--712.
%
%\bibitem{Ha2} J.~Hausen: Geometric Invariant Theory based on Weil divisors.
%Compositio Math. 140,  1518-1536 (2004)
%
\bibitem{HK} Y.~Hu, S.~Keel: Mori dream spaces and GIT.
Michigan Math. J. 48, 331--348 (2000)
%
%\bibitem{Kl}  S. L. Kleiman:  Toward a numerical theory of ampleness.
%Ann. of Math.~(2) 84, 293--344 (1966)
%
\bibitem{KP} H.~Kraft, V.L.~Popov: Semisimple group actions on the three dimensional 
affine space are linear. Comment. Math. Helvetici 60, 466--479 (1985)
%
\bibitem{KKV} F. Knop, H. Kraft, T. Vust: The Picard group of a $G$-variety,
in: Algebraische Transformationsgruppen und Invariantentheorie, DMV Seminar,
Vol. 13, Birkh\"auser, Basel (1989)
%
\bibitem{Lu}
D. Luna: Vari\'{e}t\'{e}s sph\'{e}riques de type $A$. 
Publ. Math. Inst. Hautes \'{E}tudes Sci.  No. 94 (2001), 
161--226.
%
\bibitem{LuVu}
D. Luna, Th. Vust: Plongements d'espaces homog\`{e}nes. 
Comment. Math. Helv.  58  (1983),  no. 2, 186--245.
%
\bibitem{Mu}
D.~Mumford, J.~Fogarty, F.~Kirwan:
Geometric invariant theory. Third edition. 
Ergebnisse der Mathematik und ihrer Grenzgebiete, 
34. Springer-Verlag, Berlin, 1994
%
\bibitem{Kn}
F.~Knop:
\"Uber Hilberts vierzehntes Problem 
f\"ur Variet\"aten mit Kompliziertheit eins. 
Math. Z.  213  (1993),  no. 1, 33--36.
%
\bibitem{OV} A.L.~Onishchik, E.B.~Vinberg: Lie Groups and Algebraic Groups.
Springer-Verlag, Berlin Heidelberg, 1990
%
\bibitem{Re}
N.~Ressayre: The GIT-equivalence for $G$-line bundles.  
Geom. Dedicata  81  (2000),  no. 1-3, 295--324.
%
\bibitem{Se} C.S.~Seshadri: Quotient spaces modulo reductive algebraic
groups. Ann. Math. (2) 95, 511--556 (1972)
%
\bibitem{Th}
M.~Thaddeus: Geometric invariant theory and flips.  
J. Amer. Math. Soc.  9  (1996),  no. 3, 691--723.
%
\bibitem{Ti}
D.A. Timash\"ev: Classification of $G$-manifolds of complexity $1$. 
Izv. Ross. Akad. Nauk Ser. Mat.  61  (1997),  no. 2, 127--162;  
translation in  Izv. Math.  61  (1997),  no. 2, 363--397
%
\bibitem{Vi} E.B.~Vinberg: Invariant linear connections in a homogeneous
space (Russian). Tr. Mosk. Mat. O-va 9, 191--210 (1960)
%
\bibitem{Sw1}
J.~\'Swi\c{e}cicka: Quotients of toric 
varieties by actions of subtori.  
Colloq. Math.  82  (1999),  no. 1, 105--116.
%
\bibitem{Sw2}
J.~\'Swi\c{e}cicka: 
A combinatorial construction of sets with good quotients 
by an action of a reductive group.  
Colloq. Math.  87  (2001),  no. 1, 85--102. 
%
\bibitem{Wl}
J.~W\l odarczyk:
Embeddings in toric varieties and prevarieties.
J. Algebraic Geom. 2 (1993), no. 4, 705--726.
\end{thebibliography}
\end{document}